\documentclass[times, doublespace]{article}
\usepackage{amssymb, amsmath, amsthm, mathtools, ushort, authblk}
\usepackage[all]{xy}

\newcommand{\on}{\operatorname}
\newcommand{\mc}{\mathcal}
\newcommand{\mf}{\mathfrak}
\newcommand{\mb}{\mathbf}

\newcommand{\be}{\begin{equation}}
\newcommand{\ee}{\end{equation}}

\def\lra{{{\,\,\longrightarrow\,\,}}} 
\def\lla{{{\,\,\longleftarrow\,\,}}} 
\newcommand{\ul}{\underline}
\newcommand{\ol}{\overline}
\newcommand{\ind}{\varinjlim}
\newcommand{\pro}{\varprojlim}
\newcommand{\lb}{( \hskip -.7mm (}
\newcommand{\rb}{) \hskip -.7mm )}
\newcommand{\lsb}{[\hskip -.5mm [}
\newcommand{\rsb}{]\hskip -.5mm ]}
\newcommand{\sqr}{{\sqrt{\,\,\,\,}}}

\newcommand{\Algb}{\mb {Alg}}
\newcommand{\SSchb}{\mb{Ssch}}
\newcommand{\Sspb}{\mb{Ssp}}
\newcommand{\Lrsb}{\mb{Lrs}}
\newcommand{\Schb}{\mb{Sch}}
\newcommand{\Spec}{\on{Spec}}
\newcommand{\Fschb}{\mb{Fsch}}
\newcommand{\Setb}{\mb{Set}}
\newcommand{\Ischb}{\mb{Isch}}
\newcommand{\Topb}{\mb{Top}}
\newcommand{\Sab}{\mb{SAb}}
\newcommand{\Fsetb}{{\mb{Fset}^+}}
\newcommand{\Fsetall}{{\mb {Fset}}}
\newcommand{\SVectb}{{\mb{SVect}}}
\newcommand{\Affb}{{\mb{Aff}}}

\newcommand{\even}{{\ol{0}}}
\newcommand{\odd}{{\ol{1}}}
\newcommand{\Hom}{\on{Hom}}
\newcommand{\Aut}{\on{Aut}}
\newcommand{\Der}{\on{Der}}
\newcommand{\GL}{\on{GL}}
\newcommand{\SL}{\on{SL}}
\newcommand{\cl}{{\on{cl}}}
\newcommand{\Ker}{\on{Ker}}
\newcommand{\Coker}{\on{Coker}}
\newcommand{\Tot}{\on{Tot}}
\newcommand{\don}{d}

\newcommand{\Id}{{\on{Id}}}
\newcommand{\red}{{\on{red}}}
\newcommand{\ev}{{\on{ev}}}
\newcommand{\Res}{{\on{Res}}}

\newcommand{\CC}{{\mathbb C}}

\newcommand{\GG}{{\mathbb G}}

\newcommand{\ZZ}{{\mathbb Z}}

\newcommand{\Ac}{{\mc A}}

\newcommand{\Dc}{{\mc D}} 
 
\newcommand{\Fc}{{\mc F}}

\newcommand{\Ic}{{\mc I}} 
\newcommand{\Kc}{{\mc K}} 
\newcommand{\Lc}{{\mc L}}

\newcommand{\Oc}{{\mc O}}

\newcommand{\Sc}{{\mc S}}

\newcommand{\gen}{{\mf g}}

\newcommand{\oen}{{\mf o}} 
\newcommand{\pen}{{\mf p}}

\newcommand{\ten}{{\mf t}}

\newcommand{\Cb}{{\mb C}}

\newcommand{\Eb}{{\mb E}}

\def\g{{\gamma}} 
\def\a{{\alpha}}

\newcommand{\eps}{{\varepsilon}}

\newtheorem{thm}[equation]{Theorem}
\newtheorem{cor}[equation]{Corollary}
\newtheorem{lem}[equation]{Lemma}
\newtheorem{prop}[equation]{Proposition}

\newtheoremstyle{example}{\topsep}{\topsep}%
     {}
     {}
     {\bfseries}
     {.}
     {2pt}
     {\thmname{#1}\thmnumber{ #2}\thmnote{ #3}}

   \theoremstyle{example}
   
   \newtheorem{Defi}[equation]{Definition}
   \newtheorem{rem}[equation]{Remark}
   \newtheorem{rems}[equation]{Remarks}
   
   \newtheorem{exas}[equation]{Examples}
   
  \numberwithin{equation}{subsection}

\newtheorem{exa}[equation]{Example}

\author{Mikhail Kapranov 
\thanks{ Department of Mathematics, Yale University, New Haven CT 06520,
email: mikhail.kapranov@yale.edu}
\and  Eric Vasserot
\thanks{D\'epartement de Math\'ematiques,Universit\'e Paris 7, 
175 rue de Chevaleret, 75013 Paris, France, email: vasserot@math.jussieu.fr}}

\title{Supersymmetry and the formal loop space.}

\begin{document}
\maketitle

\begin{abstract}
 For any algebraic super-manifold $M$ we define the super-ind-scheme $\Lc M$
  of formal loops and study the transgression map (Radon transform) on differential forms in this context. Applying this to the super-manifold $M=\Sc X$, the spectrum of the de Rham complex of a manifold $X$, we obtain, in particular, that the transgression map for $X$ is a quasi-isomorphism between the $[2,3)$-truncated de Rham complex of $X$ and the additive part of the $[1,2)$-truncated de Rham complex of $\Lc X$. The proof uses the super-manifold $\Sc\Sc X$ and the action of the Lie superalgebra $\mf{sl}(1|2)$ on this manifold. This quasi-isomorphism result provides a crucial step in the classification of sheaves of chiral differential operators in terms of 
  the geometry of the formal loop space.

\end{abstract}

 \tableofcontents
 
 
 \addtocounter{section}{-1}
 
 \section{ Introduction.}

 Let $X$ be a smooth complex algebraic variety. The {\em de Rham spectrum}
 $\Sc X=\Spec(\Omega^\bullet_X)$, is a super-manifold which can be seen
 as the configuration space of a supersymmetric particle moving in $X$,
 see, e.g.,  \cite{Witten}.
 The particle itself can be understood as the super-manifold
 $\mathbb A^{0|1}=\Spec\,\Lambda[\eta]$. It was pointed out by M. Kontsevich
 that the de Rham differential comes from the internal symmetry of the particle,
 i.e., from the action of the super-group  
 of automorphisms of $\mathbb A^{0|1}$.
 In fact, representations of this super-group are the same as cochain
 complexes, see \cite{Kontsevich-deformations} and Subsection
 \ref{N=1-supersymmetry-subection} below.  
 
 \vskip .2cm
 
 The functor $\Sc$ can be applied to any super-manifold, in particular, we can
 form $\Sc\Sc X=\Spec(\Omega_{\Sc X}^\bullet)$. It has a similar
 interpretation to the above,  but in terms of $\mathbb A^{0|2}=\Spec\,\Lambda[\eta_1,
 \eta_2]$ which  can be seen as  an  ``$N=2$ supersymmetric particle",   moving in $X$. 
 Mathematically, the most immediate part of $N=2$ supersymmetry  
 is the super-group $\ul{\Aut}(\mathbb A^{0|2})$, acting on $\Sc\Sc X$. Its Lie
 algebra includes the two natural differentials on $\Omega_{\Sc X}^\bullet$. 
 A remarkable feature of the $N=2$ case, lost if we pass to $\mathbb A^{0|N}$
 for $N>2$, is that $\ul{\Aut}(\mathbb A^{0|2})$ is  isomorphic to 
 the special linear super-group $SL_{1|2}$.
 Therefore $SL_{1|2}$ acts on the double complex
 $\Omega_{\Sc X}^\bullet$, and this
 action gives a detailed information about the cohomology of the rows and columns. 
 
 \vskip .2cm
 
 The goal of the present paper is to apply these ideas to the study of $\Lc X$, the
 ind-scheme of formal loops in $X$, introduced by us in \cite{KV1}. 
 We showed that $\Lc X$ possesses a nonlinear analog of a vertex algebra structure,
 called the structure of a {\em factorization semigroup}. Therefore, natural
 linear objects on $\Lc X$ give rise to vertex algebras. In particular, we showed how to obtain $\Omega^{\on{ch}}_X$,
 the chiral de Rham complex of $X$, see \cite{MSV}, from  geometry of $\Lc X$. 
 Our point of view suggests a similar interpretation of the sheaves of chiral
 differential operators (CDO) on $X$ studied in \cite{GMS1, GMS2}: a CDO can be obtained
 from an object of the determinantal gerbe of $\Lc X$ which is {\em factorizing}
 (compatible with the factorization semigroup structure). 
 In fact, the factorization structure on $\Lc X$ leads naturally
 to the factorization conditions for all sorts of geometric objects on $\Lc X$:
 functions, forms, line bundles, gerbes, etc. For functions and forms
 the factorization is understood in the additive sense, so we will speak about
 {\em additive forms} on $\Lc X$, see Definition \ref{additive-forms-definition}. 
 In \cite{GMS2},  the CDO's were classified in terms of the complex $\Omega^{[2,3)}_X$
 which is the second of the two truncated de Rham complexes below
 (on the Zariski topology of $X$):
 $$\Omega_X^{[1,2)}=\bigl\{\Omega^1_X\buildrel d\over\lra \Omega^{2, \cl}_X\bigr\},
\quad
\Omega_X^{[2,3)}=\bigl\{ \Omega^2_X\buildrel d\over\lra \Omega^{3, \cl}_X\bigr\}.$$
Here $\Omega^{i, \cl}_X$ is the sheaf of closed differential $i$-forms on $X$.
Note that $\Omega^{[1,2)}_X$
governs rings of {\em twisted differential operators} on $X$, see
 \cite{Beilinson-Bernstein}; such rings
form a stack of Picard categories,  of which the gerbe of  CDO is often said to be
 a ``higher" analog. 
On the other hand, our approach with the determinantal gerbe also leads to the
complex $\Omega^{[1,2)}$, but on $\Lc X$. 

\vskip .2cm

Our  main result, Corollary \ref{omega-2-3-corollary},
 says that $\Omega^{[2,3)}_X$ is quasi-isomorphic
 (although not isomorphic)
to a subcomplex in $\Omega^{[1,2)}_{\Lc X}$ consisting of additive forms.
The quasi-isomorphism is given by the transgression (Radon transform)
$$\tau: \Omega^p (X)\,\lra\,\Omega^{p-1} (\Lc X).$$
In fact, we prove a general statement about the full
de Rham complexes of $X$ and $\Lc X$ (Theorem 
\ref{additive-forms-transgression-theorem}), of which
Corollary \ref{omega-2-3-corollary} is a
consequence. The proof uses $N=2$ supersymmetry. 

 \vskip .2cm
 
 In \cite{KV3} we proved that additive functions $f$ on on $\Lc X$ are
identified with closed 2-forms $\omega$ on $X$ via a version of the
{\em symplectic action functional}
$$\omega\,\, \mapsto S(\omega) = d^{-1}(\tau(\omega)).$$
The argument in \cite{KV3} used vertex algebras and the result of \cite{MSV} on realization
of $\Omega^{2,\cl}_X$ as the sheaf of vertex automorphisms of the
chiral de Rham complex. Here we give a direct proof of this fact from first
principles (Theorem \ref{additive-2-forms-theorem}), by expanding $f$
around constant loops. Here $X$ can be any super-manifold. Now, viewing
additive forms on $\Lc X$ as additive functions on $\Sc\Lc X=\Lc\Sc X$, we identify
their space with $\Omega^{2,\cl}_{\Omega^\bullet_X}$ which, by the $N=2$ supersymmetry
analysis, is quasi-isomorpic to the cohomological truncation
$$\Omega^{2,\cl}_X\lra\Omega^2_X\lra\Omega^3_X\lra ...,$$
giving Theorem 
\ref{additive-forms-transgression-theorem}. 
 Note that $\Lc\Sc X$ can be seen as the space of ``super-loops"
in $X$, i.e., of maps from a super-thickening of the punctured formal disk. 
In Subsection \ref{super-curves-factorization-subsection} we show how
$\Lc\Sc^NX$  gives rise to a factorization semigroup on any $(1|N)$-dimensional
super-curve.

\vskip .2cm

This paper was originally intended as an appendix to \cite{KV4}
but it seemed better to us to write it separately,
collecting together the aspects of the theory 
with less emphasis on categorical issues. 
These categorical issues, i.e.,  complete yoga of 
  factorization as applied to not just functions and forms  on $\Lc X$, but
  $\Dc$-modules, line bundles,  gerbes etc.,
 form  the natural subject of  \cite{KV4}, whose place in the logical
 order is after the present paper. In fact, factorizing
 gerbes can be given a de Rham-type description,
much in the spirit of the book \cite{Brylinski} by Brylinski. This description leads 
to a direct identification of the gerbe of CDO with the gerbe corresponding to
the additive part of the complex $\Omega^{[1,2)}_{\Lc X}$. The results of
the present paper, identifying this additive part with $\Omega^{[2,3)}_X$,
provide then a clear explanation of the classification of \cite{GMS2}
from the first principles.

\vskip .2cm

Here is a brief outline of the paper. In Section 1 we provide 
the necessary background for for scheme-theoretic
algebraic geometry in the super-setting. As there seems to be
 no systematic reference in the required generality,
 we had to give a somewhat longer treatment.
 Section 2 is devoted to the discussion of extended supersymmetry
 from the point of view of super-version of schemes of infinitesimally
 near points in the spirit of A. Weil \cite{Weil}. Here we analyze representations
 of $\ul{\Aut}(\mathbb A^{0|2})=SL_{1|2}$ as double complexes
 with appropriate partial contracting homotopies. 
 In Section 3, we discuss 
 the formalism of super-ind-schemes, quite parallel to that of usual
 ind-schemes. 
 In Section 4 we define formal loop spaces in the
 super-setting while in Section 5 we discuss their factorization
 structure. The formalism of factorization data which we discuss
 differs slightly from that of \cite{BD1}; it is  
 better adapted to
 studying coherence conditions needed for factorizing line bundles, gerbes
 etc. Finally, in Section 6 we prove our main results: first about
 additive functions, then about additive forms.

 \vskip .2cm
 
 We are grateful to D. Osipov and A. Zheglov for pointing out some inaccuracies
 in \cite{KV1}. We correct these inaccuracies in the present paper. We would also like to
 thank D. Leites and A. Zeitlin for pointing out some classical references
 dealing with supersymmetry. The first author acknowledges the support of an NSF grant
 and of the Universit\'e Paris-7, where a part of this work was written. 
 
 
 \section{Superschemes.} 
 
 \subsection{Basic definitions.}

 We start by discussing basic concepts of algebraic geometry
 in the super situation, following \cite{Leites-spectra, Manin,
 Voronov-Manin-Penkov}. See also \cite{Deligne-Morgan}
 for a general background in a more differential-geometric context.
 
 First of all, recall that a {\em ringed space}
  is a pair $X=(\underline{X},\Oc_X)$ where $\underline{X}$ is a topological
  space,
and $\Oc_X$ is a sheaf of rings, not necessarily commutative,
on $\underline{X}$. A morphism 
$f: X=(\underline{X}, \Oc_X)\to Y=(\underline{Y}, \Oc_Y)$ of ringed spaces
consists of  a continuous map of spaces $f_\sharp: \underline{X}\to\underline{Y}$,
and a morphism of sheaves of rings 
$f^\flat: f^{-1}_\sharp(\Oc_Y)\to\Oc_X$
on $\underline{X}$.

An {\em open embedding} of ringed spaces is a morphism $f$ such that
$f_\sharp$ is an open embedding of topological spaces, while
$f^\flat$ is an isomorphism of sheaves of rings. 

A {\em locally ringed space} is a ringed space $X=(\underline{X},\Oc_X)$
such that each stalk $\Oc_{X,x}, x\in \underline{X}$, is a local ring.
A {\em local morphism} of locally ringed spaces is a morphism
$f$ as above such that each morphism of stalks
$f^\flat_x: \Oc_{Y, f_\sharp(x)}\to \Oc_{X,x}$ is a local
homomorphism of local rings, i.e., takes the maximal ideal of one ring
into the maximal ideal of the other. For example, an open
embedding of locally ringed spaces is always a local morphism. 
We denote by $\Lrsb$ the category
of locally ringed spaces and their local morphisms. 

We also denote by $\Schb$ the category of schemes. Recall that
$\Schb$ is a full subcategory in $\Lrsb$. In particular, for any commutative
ring $R$ we have the scheme $\Spec(R)$ whose underlying topological
space (i.e., the set of prime ideals in $R$ with the Zariski topology)
will be denoted $\ul{\Spec}(R)$. 

\vskip .2cm

We denote by $\Sab$ the symmetric monoidal category of $\ZZ/2$-graded abelian groups
$A=A_\even \oplus A_\odd$. The symmetry transformation $A\otimes B\to B\otimes A$
in this category is given by the Koszul sign rule: 
$$a\otimes b\mapsto (-1)^{\deg(a)\deg(b)}b\otimes a$$
on homogeneous elements. We denote by 
 $\Pi: \Sab\to\Sab$  the functor of
 change of parity: $(\Pi A)_\even = A_\odd$ and vice versa.

 Recall that a {\em super-commutative ring} is a commutative ring object
 in the symmetric monoidal category $\Sab$. Explicitly, it is
 a $\ZZ/2$-graded ring $R=R_\even \oplus R_\odd$
such that $ab=(-1)^{\deg(a)\deg(b)}ba$ for homogeneous elements. 
 The following is then clear.

\begin{prop} Let $R$ be a super-commutative ring. Then:

(a) A $\ZZ/2$-graded ideal $\mf p=\mf p_\even\oplus\mf p_\odd\subset R$
is prime (in the sense that $R/\mf p$ has no zero-divisors), if and only if
$\mf p_\even$ is a prime ideal in $R_\even$. In this case
$\mf p_\odd=R_\odd$. 

(b) The Jacobson radical of $R$ is equal to the sum $\sqrt{R_\even}\oplus R_\odd$.

(c)  $R$ is local with maximal ideal $\mf m$ if and only if $R_\even$ is local
with maximal ideal $\mf m_\even=\mf m\cap R_\even$.  
\qed
\end{prop}

A {\em super-space} is a locally ringed space $(\ul X,\Oc_X)$ where $\Oc_X$ is equiped
with a $\ZZ/2$-grading $\Oc_X=\Oc_{X,\even}\oplus\Oc_{X,\odd}$
making it into a sheaf of super-commutative rings.
A {\em morphism of super-spaces} is a local morphism of locally ringed spaces
$f=(f_\sharp, f^\flat)$ such that $f^\flat$ preserves the $\ZZ/2$-grading.
We denote by $\Sspb$ the category of super-spaces.

A super-space $X= (\ul X,\Oc_X)$ is called a {\em super-scheme} if $(\ul X,\Oc_{X,\even})$ is a 
scheme and $\Oc_{X,\odd}$ is a quasi-coherent of $\Oc_{X,\even}$-module.
We denote by $\SSchb\subset \Sspb$ the full subcategory formed
by super-schemes.

\vskip .2cm

In particular for any super-commutative ring 
$R$ we have a super-scheme $\Spec(R)$.  Its underlying space is
$\ul{\Spec}(R)=\ul{\Spec}(R_\even)$, with
$\Oc_{\Spec(R)}=\Oc_{\Spec(R),\even}\oplus \Oc_{\Spec(R),\odd}$
where $\Oc_{\Spec(R), \even}$ is the structure sheaf of $\Spec(R_\even)$,
while $\Oc_{\Spec(R), \odd}$ is the quasi-coherent sheaf of 
$\Oc_{\Spec(R), \even}$-modules corresponding to the $R_\even$-module $R_\odd$. 
Super-schemes of the form $\Spec(R)$ will be called {\em affine}. 
It is clear that every super-scheme is locally isomorphic to
an affine super-scheme.

\vskip .2cm

Given a super-scheme $X$, a {\em quasi-coherent sheaf} of $\Oc_X$-modules
is a $\ZZ/2$-graded sheaf $\Fc=\Fc_\even\oplus\Fc_\odd$ 
which is quasi-coherent as a sheaf of $\Oc_{X,\even}$-modules. 
For any super-commutative ring $R$ quasi-coherent sheaves
on $\Spec(R)$ are in bijection with $\ZZ/2$-graded $R$-modules.

For a morphism of super-commutative algebras $A\to B$ we denote by
$\Omega^1(B/A)$ the $\ZZ/2$- module 
of K\"ahler differentials
of $B$ over $A$ understood in the super-sense,
so that $d: B\to\Omega^1(B/A)$ preserves the $\ZZ/2$-grading, annihilates the
image of $A$, and
satisfies the super-Leibniz rule. Alternatively,
\be\label{omega-1-kernel-equation}
\Omega^1 ({B/A})=I/I^2, \quad I\,=\,\Ker\bigl\{ B\otimes_A B\buildrel
\on{mult.}\over\lra B\bigr\}.
\ee
For a morphism of super-schemes $X\to Y$ we have then
the quasi-coherent sheaf $\Omega^1_{X/Y}$ on $X$.

\vskip .2cm

Given a super-scheme $X$ and a quasi-coherent sheaf $\Ac$ of
super-commutative $\Oc_X$-algebras, we have a super-scheme
$\Spec_X(\Ac)\to X$ obtained by gluing affine schemes
$\Spec(\Ac(U))$ for open affine sub-schemes $U\subset X$. 
A morphism $Y\to X$ of super-schemes is called {\em affine},
if it is isomorphic to one of the form $\Spec_X(\Ac)\to X$.
Note the particular case when $\Ac = \Oc_X/\Ic$
is the quotient of $\Oc_X$ by a sheaf of ideals. In this case
$Y=\Spec_X(\Ac)$ is called a {\em closed sub-super-scheme} in $X$,
and any morphism isomorphic to $Y\to X$ of this type is called a
 {\em closed embedding} of super-schemes. 
 An {\em immersion} is a morphism of super-schemes which can be represented as
the composition of an open embedding followed by a closed embedding. 

\vskip .2cm

A super-scheme $X$ will be called {\em quasi-compact}, if the topological
space $\ul{X}$ is quasi-compact, i.e.,  each open covering of $\ul{X}$
has a finite sub-covering.
For example every affine super-scheme is quasi-compact. 

\vskip .2cm

As in the case of ordinary schemes, we have the following fact.

\begin{prop}

(a) The category $\SSchb$ has finite projective limits,  in particular, finite products 
and
fiber products.

(b) Let $I$ be a filtering poset and  $(X_i)_{i\in I}$ be a
 projective system of super-schemes with structure morphisms 
  $u_{ij}: X_j\to X_i$, given for $i\leqslant j$. If all $u_{ij}$ are affine morphisms,
  then  the limit $\pro_{i\in I}^{\SSchb} X_i$ exists. Denoting this limit by $X$,
  we have that the natural projection $p_i: X\to X_i$ is affine for any $i$,
  in fact
$$X\,\,=\,\, \Spec_{X_i}\biggl( \ind_{j\geqslant i}\, u_{ij*}\Oc_{X_j}\biggr).$$ 
  Moreover, we have $\ul{X}=\pro_{i\in I}^{\Topb}\ul{X}_i$. 
  \qed
\end{prop}

\noindent {\sl Proof:} (a) In any category, existence of finite projective limits
is equivalent to the existence of finite products and fiber products. Now, for
affine super-schemes, the fiber product of
$$\Spec(A)\lra \Spec(C)\longleftarrow \Spec(B)$$
is found as $\Spec(A\otimes_C B)$, like for ordinary schemes. After that,
 fiber products of arbitrary super-schemes are defined by gluing affine charts of the
kind described. 

\vskip .2cm

(b) The argument is identical to \cite{EGA-IV-8-15} (8.2.3) (existence of the limit and
its realization as a relative spectrum) and (8.2.10) (description of $\ul{X}$). \qed

\vskip .3cm

Given a super-scheme $X = (\ul{X} ,\Oc_X)$ its {\em even part} is defined to be the scheme
\be X_{\on{even}}\,\,=\,\, \bigl(\ul{X}, \Oc_{X,\even}/(\Oc_{X,\odd}^2)\bigr),
\ee
while the corresponding {\em reduced scheme} is 
\be\label{reduced-scheme-equation}
X_\red\,\,=\,\,\left(\ul{X}, \Oc_{X,\even}\bigr/\sqrt{\Oc_{X,\even}}\right),
 \ee
 similarly to the case of ordinary schemes. 
 
 \vskip .2cm
 
 From now on we work over the field $\CC$ of complex numbers.
All rings will be assumed to contain 
$\CC$ and all super-schemes will be super-schemes over $\CC$. 
 We denote by $\SVectb$ the symmetric monoidal category of
 $\ZZ/2$-graded $\CC$-vector spaces, and by
  $\Algb$ the category of super-commutative $\CC$-algebras. 
We also denote  $\Affb$ the category of affine super-schemes,
 i.e., the dual category of $\Algb$.

 An affine super-scheme (over $\CC$) is said to be of {\em  finite type} if it is
isomorphic to $\Spec(R)$ where $R$ is a finitely generated super-commutative
$\CC$-algebra. More generally a {\em super-scheme of finite type} 
is a super-scheme  which
 can be covered by finitely many affine super-schemes of 
finite type. Let 
 $\Fschb$ 
be the full category of super-schemes of finite type.

\begin{exas}\label{affine-space-examples}
 (a) For $d_1, d_2\geqslant 0$ we denote by $\CC^{d_1|d_2}\in\SVectb$ the
 coordinate $\ZZ/2$-graded space with $d_1$ even dimensions and $d_2$ odd
 dimensions. For $R\in\Algb$ we denote by $R^{d_1|d_2}$ the $\ZZ/2$-graded
 $R$-module $R\otimes \CC^{d_1|d_2}$. If $R$ is local, then any
 finitely generated projective $R$-module $M$ is free, i.e.,  isomorphic to $R^{d_1|d_2}$
 for a unique pair $(d_1|d_2)$ which is called the {\em rank} of $M$.
 If $R$ is finitely generated, then a finitely generated projective $R$-module
 is free, locally on the Zariski topology of $\ul{\Spec}(R)$, so its rank
 is a locally constant function $\ul{\Spec}(R)\to\ZZ_+\times\ZZ_+$.

 We define $\mathbb A^{d_1|d_2}$, the {\em affine super-space} of dimension
$d_1|d_2$ to be the super-scheme
$$\mathbb A^{d_1|d_2}=\Spec \bigl(\CC[x_1,\dots,x_{d_1}]\otimes
\Lambda[\xi_1,\dots,\xi_{d_2}]\bigr).$$

\noindent (b) It will be convenient to use the following unified notation.
 Given  $N=d_1+d_2$ generators $a_1, ..., a_N$ of which $d_1$ are even and $d_2$ are odd, 
 we will simply write $\CC[a_1, ..., a_N]$
for the  tensor product of the polynomial algebra on the even generators and the exterior
algebra on the odd generators. We also write $\CC\lsb a_1, ..., a_N\rsb$ for the
completion of $\CC[a_1, ..., a_N]$ with respect to the ideal $(a_1, ..., a_N)$,
which is the tensor product of the formal power series algebra on the even generators
and the exterior algebra on the odd generators. 

\vskip .2cm

\noindent (c) For any $d_1,d_2\geqslant 0$ we have the group super-scheme $GL_{d_1|d_2}$ such that
for a super-commutative algebra $R$, the group $GL_{d_1|d_2}(R)$ consists of $\ZZ/2$-homogeneous
automorphisms of the $R$-module $R\otimes \CC^{d_1|d_2}$.
Such automorphisms can be represented by block matrices over $R$ of format $(d_1+d_2)\times(d_1+d_2)$
$$g=\begin{pmatrix} A&B\\C&D\end{pmatrix}
$$
with entries of $A,D$ belonging to $R_\even$, entries of $B,C$ to $R_\odd$,
and $A,D$  invertible. 
The group super-scheme $GL_{d_1|d_2}$ is called the {\em general linear group} of format $d_1|d_2$. In particular, $GL_{1|0}=\GG_m$
is the multiplicative group. The {\em Berezin determinant}
 is a morphism of group super-schemes
$\on{ber}: GL_{d_1|d_2}\to\GG_m$ which  on $R$-points sends a matrix $g$ as above to
$$ \on{ber}(g) = \det(A-BD^{-1}C) /\det(D) \,\,\in\,\, R^*_\even = \GG_m(R).$$ 
Its kernel is denoted by $SL_{d_1|d_2}$ and called
the {\em special linear group} of format $d_1|d_2$. 

\end{exas}


\subsection{Smooth and \'etale morphisms.}

Let $f: X\to Y$ be a morphism of super-schemes. As in the classical (even) case,
we say that $f$ is {\em locally of finite presentation}, if $\Oc_X$ is, locally
on the Zariski topology of $X$, finitely presented as an $\Oc_Y$-algebra,
i.e., given by finitely many generators and relations. 

A morphism of super-commutative algebras $u:R\to R'$ is called a {\em simple extension},
if $u$ is surjective, and $I=\Ker(u)$ satisfies $I^2=0$. 
Recall that each super-scheme $X$ gives a covariant functor
$h^X: \Algb\to\Setb$, sending $R$ to $\Hom_{\SSchb}(\Spec(R), X)$. 

\begin{Defi}\label{formal-smoothness-definition}
 (a) Let $f: h\to h'$ be a morphism of covariant functors $\Algb\to\Setb$.
We say that $f$ is {\em formally smooth} (resp. {\em formally \'etale}) if, for any
simple extension $R\to R'$ the natural morphism
$$h(A)\,\lra\, h'(A)\times_{h'(A')} h(A')$$
is surjective (resp. bijective). 

(b) A morphism $f: X\to Y$ of super-schemes is called {\em formally smooth}
(resp. {\em formally \'etale}), if the corresponding morphism of  functors
$h^X\to h^Y$ is formally smooth (resp. formally \'etale). 
\end{Defi} 

A morphism of super-schemes is called {\em smooth}, if
it is formally smooth and locally of finite presentation.  

\begin{prop} (a) Let $\psi: A\to B$ be a morphism of super-commutative algebras
such that $\psi^*: \Spec(B)\to\Spec(A)$ is a formally smooth morphism of
super-schemes. Then the $B$-module $\Omega^1(B/A)$ is projective.

(b) Let $f: X\to Y$ be a smooth morphism of super-schemes.
Then $\Omega^1_{X/Y}$ is locally free, as a sheaf of $\Oc_X$-modules.
\end{prop} 

In particular,   
the rank of $\Omega^1_{X/Y}$ is a locally constant function on $\ul{X}$
with values in $\ZZ_+\times\ZZ_+$
denoted by $\dim(X/Y)$ and called the {\em relative dimension} of $X$ over $Y$.
If $f$ is \'etale, then $\dim(X/Y)$ is identically equal to 0. 

\vskip .2cm

\noindent {\sl Proof of the proposition:} (a) The classical argument
(contained in a more general form in \cite{EGA-0-14-23}
(19.5.4.1)), is completely
formal and goes in our case as follows.

Any $\ZZ/2$-graded $B$-module $Q$ gives a super-commutative $B$-algebra
$B\oplus Q$ with $Q^2=0$ and with the multiplication of $B$ and $Q$ given by
the module structure. A $B$-module homomorphism $u: \Omega^1(B/A)\to Q$
is the same as a $Q$-valued derivation $\delta: B\to Q$ vanishing on $A$,
and this gives a homomorphism of $B$-algebras
$$(\Id,\delta): \,\,B\lra B\oplus Q, \quad b\mapsto (b,\delta(b)).$$
To prove that $\Omega^1(B/A)$ is projective, let $s: P\to Q$ be a surjective
morphism of $B$-modules, and $u: \Omega^1(B/A)\to Q$ be any morphism of $B$-modules.
We prove that $u$ can be lifted to a $v: \Omega^1(B/A)\to P$.
Indeed, $(\Id\oplus s): B\oplus P\to B\oplus Q$ is a simple extension of super-commutative
algebras, and we have a commutative square
$$\begin{matrix}
&A&{\buildrel(\psi,0)\over\lra}&B\oplus P&\cr
\psi&\big\downarrow&&\big\downarrow&\Id\oplus s\cr
&B&{\buildrel (\Id,\delta)\over\lra}&
B\oplus Q
\end{matrix}$$
So by the condition that $\psi^*: \Spec(B)\to\Spec(A)$ is formally smooth,
we find that there is an algebra homomorphism $w: B\to B\oplus P$
splitting the square into two commutative triangles. The second component
of $w$ gives a derivation $B\to P$ lifting $\delta$, i.e., a homomorphism
$v: \Omega^1(B/A)\to P$ as claimed.

\vskip .2cm

(b) This follows from the fact (proved in  the same way as in the
commutative case) that a finitely presented projective
module over any super-commutative algebra is locally free. \qed

\vskip .3cm

A {\em smooth algebraic super-variety} is a super-scheme $X$ of finite type (over $\CC$)
such that $X\to\Spec(\CC)$ is a smooth morphism.
In this case $X_{\on{even}}$
is a smooth algebraic variety over $\CC$ in the usual sense.
We write $\dim(X)$ for $\dim(X/\CC)$.
If $X$ is irreducible (i.e., $\ul{X}$ is an irreducible topological space),
then this function is constant, so $\dim(X)=(d_1|d_2)$ for some $d_1, d_2\in\ZZ_+$,
and  $\dim(X_{\on{even}})=d_1$. If $\dim(X)=(1|N)$, we say that
$X$ is a {\em super-curve}.

\begin{prop} Let $X\to Y$ be a smooth morphism, $x\in X(\CC)$ be such that
$\dim(X/Y)=(d_1|d_2)$ at $x$. Then there are Zariski open sets $U\subset X$
containing $x$  such that  
there is a morphism of $Y$-schemes $\phi: U\to Y\times \mathbb A^{d_1|d_2}$
which is \'etale. 
 \end{prop}
 
 \noindent {\sl Proof:} analogous to the purely even case which
 is proved in  \cite{EGA-IV-16-21} (17.11.4).
 As $f$ is locally of finite presentation, we find $U'$ and
 $V'$ with $x\in U'$, $f(x)\in V'$,  $f(U')\subset V'$ so that there a morphism of $V'$-schemes
 $i: U'\to V'\times \mathbb A^{D_1|D_2}$ which is a closed embedding. Then, we can choose
 a subset of the coordinates $x_i, \xi_j$ on $\mathbb A^{D_1|D_2}$ such that
 $dx_i, d\xi_j$ from that subset form a set of free generators of $\Omega^1_{U'/V'}$
 in some $U\subset U'$ containing $x$. The projection on the coordinate affine subspace
  $\mathbb A^{d_1|d_2}$ corresponding to this subset, is the  \'etale morphism
  required. \qed



 
 \section{Infinitesimally near points and supersymmetry.}
 
 \subsection {Infinitesimally near points.}

\begin{Defi}\label{inf-near-morphisms-def}
 Let $u,u': S\to X$ be morphisms of super-schemes. We say that $u$ and $u'$
are {\em infinitesimally near}, if $u=u'$ on $S_{\red}$.
\end{Defi}

In this section, we want to study super-schemes which classify such morphisms for a particular
class of super-schemes $S$. We start with general categorical remarks. 

\vskip .2cm

Let $\Cb$ be any category with finite products.
Given two objects $B,C$ of $\Cb$, we have the contravariant functor
$$\Cb\to\Setb,\quad T\mapsto\Hom_{\Cb}(T\times B,C).$$
If this functor is representable, then the representing object of $\Cb$
is denoted by $\underline\Hom(B,C)$ and is called the {\em internal Hom} from $B$ to
$C$. Note that if $B=C$, then $\underline\Hom(B,B)$ is a semigroup 
object in $\Cb$. Indeed, for every $T$, the set $\Hom_\Cb(S\times B,B)$ is a 
semigroup with unit being the canonical projection 
$S\times B\to B$. Further, consider the functor associating to $T$ the set of
invertible elements in the semigroup $\Hom_\Cb(T\times B,B)$.
If this functor is representable, then the representing object in $\Cb$ is 
denoted by $\underline\Aut(B)$ and is called the {\em internal automorphism
group} of $B$. It is a group object of $\Cb$.
 
 \vskip .2cm
 
 We now specialize to $\Cb =\SSchb$. Let $\oen$ be a finite dimensional local 
 super-commutative $\CC$-algebra.
 
  \begin{prop}\label{near-points-proposition}
(a) For any super-scheme $S$ we have an identification of
super-spaces $T\times\Spec(\oen)=(\underline T,\Oc_T\otimes\oen)$.

(b) Let $X$ be any super-scheme. Then there exists
the internal Hom superscheme $X^\oen = \underline\Hom(\Spec(\oen),X)$
  representing the functor
  $$T\mapsto \Hom\bigl(T\times\Spec(\oen),X\bigr).$$

(c) If $U\subset X$ is open then $U^\oen=X^\oen\times_XU$. In particular,
 $U^\oen$ is open in $X^\oen$. 

(d) We have $(X^{\oen_1})^{\oen_2}=X^{\oen_1\otimes\oen_2}$.

(e) The functor $X\mapsto X^\oen$ takes closed embeddings to closed embeddings.
\end{prop}

The super-scheme 
$X^\oen$  will be called the {\em superscheme of 
$\oen$-infinitesimaly near points} of $X$. This terminology and notation
is borrowed from A. Weil \cite{Weil}. 

\vskip .2cm
 
 \noindent{\sl Proof :}
 (a) Clear since $\oen$ is a finite dimensional local supercommutative algebra, and so
its maximal ideal consists of nilpotent elements.
 
 \vskip .2cm
 
 (b)  Assume that $X=\Spec(A)$ is affine.
Choose a basis $(e_i)_{i\in I}$ of homogeneous elements of
$\oen$ with the following properties. First, we assume that
 $I$ has a distinguished element $0$, and 
 $e_0=1$. Second, we assume that all $e_i, i\neq 0$, 
lie in the maximal ideal of $\oen$.  After this, write the multiplication
law in $\oen$ as
  $$e_ie_j=\sum_kc_{ij}^ke_k.$$
Define a super-commutative algebra $A^\oen$ containing $A$
generated by symbols
$a[i]$, with $a\in A$, $i\in I$, subjects to the relations
\begin{equation}\label{A^o-generators-eq}
  \begin{split}(ab)[k]=\sum_{i,j}c_{ij}^ka[i]b[j],\quad
(a+\lambda b)[i]=a[i]+\lambda(b[i]),\\ a,b\in A,\quad i\in I,\quad
\lambda\in\CC.
  \end{split}
  \end{equation}
Here the degree of $a[i]$ is the sum of the degrees of $a$ and $e_i$.
Notice that the correspondence $a\mapsto a[0]$ defines an an algebra embedding
$A\subset A^\oen$, because $(ab)[0]=a[0]\,b[0]$ for each $a,b$.
We claim that
$$\Hom_{\Algb}(A^\oen,R)=\Hom(A,R\otimes\oen),$$
for each super-commutative $\CC$-algebra $R$.
Indeed, given $f:A\to R\otimes\oen$, we expand it in the form
$$f(a)=\sum_if_i(a)\otimes e_i,\quad a\in A.$$
Then we form the map 
$$\phi:A^\oen\to R,\quad a[i]\mapsto f_i(a).$$
Note that the relations in $A^\oen$ insure that $\phi$ is a well-defined
homomorphism.
This proves  (b) for $X=\Spec(A)$, and an affine super-scheme
$T=\Spec(R)$.
This implies the equality for any $T$ in virtue of part (a), because the two
functors are sheaves on the Zariski topology of $T$.
  
  \vskip .2cm

We next prove  (c) in the particular case where $X=\Spec(A)$ is affine
and $U=\Spec(A[s^{-1}])$ is a principal open subset.
We identify the functors represented by $U^\oen$ and $X^\oen\times_X U$
on an affine super-scheme $T=\Spec(R)$.
First, $\Hom(T,U^\oen)$ consists of algebra homomorphisms
$f:A\to R\otimes\oen$ such that $f(s)$ is invertible in $R\otimes\oen$.
  
  \vskip .2cm

 Next, $\Hom(T,X^\oen\times_XU)$ consists of $f$'s as before
such that $f_0(s)$ is invertible in $R$.
They coincide by Nakayama's lemma.
Having proved (c), and (b) for an affine $X$, we deduce (b), (c) 
for any $X$ by glueing along open parts.
 
 \vskip .2cm

Part (d) is clear because the two super-schemes represent the same functor.
 
 \vskip .2cm

Finally, let us prove   (e). It is enough to observe that 
a surjective algebra homomorphism $A\to B$ yields a surjective
algebra homomorphism $A^\oen\to B^\oen$.
\qed
 
\vskip 3mm
 
 Consider the particular case where $X=\Spec(\oen)$. 
Then
  $$\Spec(\oen)^\oen=\underline\Hom(\Spec(\oen),\Spec(\oen))$$
is a semigroup super-scheme.

\begin{prop}
The object $\underline\Aut(\Spec(\oen))$ exists. It is an open
subgroup-super-scheme in the semi-group super-scheme $\Spec(\oen)^\oen$. 
\end{prop}

We abbreviate $G_\oen=\underline\Aut(\Spec(\oen))$. 
Its Lie algebra 
$\gen_\oen=\Der_\CC(\oen,\oen)$
is just the Lie super-algebra of derivations of the super-commutative 
$\CC$-algebra $\oen$.
By construction $G_\oen$ acts on $X^\oen$ for any $X$.

\vskip .3cm
  
  \noindent{\sl Proof:}
By construction the algebra $\oen^\oen$ is generated by the elements
$u_{ij}=e_i[j]$ of degree equal to the sum of the degrees
of $e_i$ and $e_j$. We have therefore a matrix $U=(u_{ij})_{i,j\in I}$
over $\oen^\oen$.
Let $\oen^\oen[U^{-1}]$ be the localization of $\oen^\oen$ obtained
by adjoining the matrix elements of $U^{-1}$. More precisely, we have
a decomposition $I=I_\even\sqcup I_\odd$ according to the parities of the $e_i$'s.
The matrix $U$ has the corresponding block decomposition
$$\begin{pmatrix} U_{\even\even}& U_{\even\odd}\cr U_{\odd\even}&U_{\odd\odd}\end{pmatrix},$$
 and elements of $U_{pq}$ have the $\ZZ/2$-degree $p+q$.
Therefore the algebra $\oen^\oen[U^{-1}]$ 
is obtained by inverting the determinants of the even matrices
$U_{\even\even}$ and $U_{\odd\odd}$. Our proposition
  is implied by the following. 
  
  \begin{lem}
The functor assigning to a given super-scheme $T$ the set of invertible
elements in $\Hom_\Schb(T\times\Spec(\oen),\Spec(\oen))$ is represented by
$\Spec(\oen^\oen[U^{-1}])$.
\end{lem}

\noindent{\sl Proof  of the lemma:}
Assume that $T=\Spec(R)$.
Then 
$$\Hom_\Schb\bigl(T\times\Spec(\oen),\Spec(\oen)\bigr)\,\,=\,\,
  \Hom_\Algb(\oen,R\otimes\oen).$$
To every homomorphism $f:\oen\to R\otimes\oen$ we associate the matrix
$(f_{ij})$ over $R$ such that
$$f(e_i)=\sum_jf_{ij}\otimes e_j.$$
Then the composition in
the semigroup
$\Hom_\Algb(\oen,R\otimes\oen)$
corresponds to the multiplication of matrices.
Next, we have an identification
$$\Hom_\Algb(\oen,R\otimes\oen)=\Hom_\Algb(\oen^\oen,R)$$
which takes $f$ to the map $u_{ij}\mapsto f_{ij}$.
Therefore $f$ is invertible if and only if the matrix $(f_{ij})$ is invertible
over $R$, which means that the matrix $U$ is mapped to an invertible matrix.
\qed

  
  \subsection {N=1 supersymmetry.}\label{N=1-supersymmetry-subection}
  
Let $\oen=\Lambda[\eta]$ be the exterior algebra in one variable, so that
$\Spec(\oen)=\mathbb A^{0|1}$. 
\vskip .2cm

Let us describe the group super-scheme 
$$G_{\Lambda[\eta]}=\underline\Aut(\mathbb A^{0|1})$$
and its Lie algebra.
For any super-commutative algebra $R$ the group
$$\Hom(\Spec(R),G_{\Lambda[\eta]})$$
consists of changes of variables of the form
$\eta\mapsto a+b\eta$, where $a\in R_\odd$ is arbitrary 
and $b\in R_\even$ is invertible.
The even part of the super-group is $\GG_m$.
The Lie superalgebra $\Der \, {\Lambda[\eta]}$ consists of the derivations
$$(a+b\eta){\don\over\don\eta}\,,\quad a,b\in\CC.$$
So its basis is formed by
$$D={\don\over{\don}\eta}\,,\quad
\Theta=\eta{\don\over\don\eta},$$
with $D$ odd and $\Theta$ even, subject to the relations
\begin{equation}
[D,D]=[\Theta,\Theta]=0,\quad [\Theta,D]=D.
\end{equation}

The following fact was pointed out by M. Kontsevich \cite{Kontsevich-deformations}.

\begin{prop}\label{kontsevich-complex-prop}

Let $V=V_\even\oplus V_\odd$ be a super-vector space. Then an action of
$G_{\Lambda[\eta]}$ on $V$ is the same as a structure of a cochain 
complex on $V$, i.e., a choice of a $\ZZ$-grading 
$V=\bigoplus_{n\in\ZZ}V^n$ such that
$$V_\even = \bigoplus_{n\in 2\ZZ}V^n,\quad V_\odd = \bigoplus_{n\in 1+2\ZZ}V^n,$$
and a differential $d:V\to V$ of degree 1 with $d^2=0$.
\end{prop}

\noindent{\sl Proof :}
The action of $\GG_m\subset G_{\Lambda[\eta]}$ gives the grading, 
so that the action of
$\Theta$ is given by $\Theta=n$ on $V^n$.
The action of
$D\in\Der\, {\Lambda[\eta]}$ gives $d$. The fact that $d$ is of degree 1 follows
from the relation $[\Theta,D]=D$. 
\qed

  \vskip .3cm

Given a super-scheme $X$, we denote
 $X^{\Lambda[\eta]}=\underline\Hom(\mathbb A^{0|1}, X)$ by $\Sc X$ and
call it the {\em De Rham spectrum}  of $X$. 
The super-scheme $\mathbb A^{0|1}$
can be called the {\em $N=1$ supersymmetric particle}
in the same sense as $\Spec(\CC)$ can be
thought as representing a point particle.  The super-scheme $\Sc X$ 
is therefore
  the configuration space of an $N=1$ supersymmetric particle moving in $X$.

Denote by $\Omega^1_X$  the sheaf of K\"ahler differentials on $X$, 
and $\Omega^\bullet_X=S^\bullet(\Pi\Omega^1_X)$ be the sheaf of 
differential forms on $X$.
Here $\Pi$ is the change of parity functor.
The derivation $\don:\Oc_X\to\Omega^1_X$ gives rise to a derivation
$\don$ on $\Omega^\bullet_X$ of degree one and square zero. 
Let $\varpi: \Sc X\to X$ be the projection.
 
\begin{prop}
We have $\varpi_*\Oc_{\Sc X}=\Omega^\bullet_X$,
with the structure of a complex on the right hand side corresponding to the
$G_{\Lambda[\eta]}$-action on the left hand side.
\end{prop}

\noindent{\sl Proof :}
Let $X=\Spec(A)$.  A basis of the algebra $\Lambda[\eta]$ consists of two elements
$e_0=1$ and $e_1=\eta$. Therefore $A^{\Lambda[\eta]}$ is the algebra
generated by $a[0]=a$, $a[1]$, given for $a\in A$  and subject to the relations
$$(ab)[1]=a(b[1])+a[1]b,\quad a,b\in A.$$
These relations are identical to those defining $\Omega^1_X$, with
$a[1]$ corresponding to $\don a$.
Further, $\deg(a[1])=\deg(a)=1$.
 So taking the super-commutative algebra $A^{\Lambda[\eta]}$
amounts to forming the symmetric algebra of $\Pi\Omega^1_X$.
\qed
 
 \begin{exa} In particular, the De Rham differential in
$\Omega^\bullet_X$ corresponds to a vector field $D$ on $\Sc X$. Assume that $X$ is a smooth
super manifiold with local coordinates $x_1,\dots,x_n$. Then on $\Sc X$
we have local coordinates $x_1,\dots,x_n,\xi_1,\dots,\xi_n$ where
$\xi=\don x_i$. The vector fields $D$ and $\Theta$ have the form
$$D=\sum\xi_i{\partial\over\partial x_i},\quad
\Theta=\sum\xi_i{\partial\over\partial\xi_i}.$$
 \end{exa}

 
 \subsection{$N=2$ supersymmetry.}
 
 Let $\oen=\Lambda[\eta_1,\eta_2]$ 
be the exterior algebra in two variables. The super-scheme
$\Spec(\oen)=\mathbb A^{0|2}$ can be called the
{\em $N=2$ supersymmetric particle}. 
The group super-scheme 
$$G_{\Lambda[\eta_1,\eta_2]}=\underline\Aut(\mathbb A^{0|2})$$
and its Lie algebra $\Der(\Lambda[\eta_1,\eta_2])$
possess remarkable symmetry properties. By definition, 
for any super-commutative algebra $R$ the group
$$\Hom(\Spec(R),G_{\Lambda[\eta_1,\eta_2]})$$
consists of change of variables of the form
$$\gathered
\eta_1\mapsto a^1+b^1_1\eta_1+b^1_2\eta_2+c^1\eta_1\eta_2\cr
\eta_2\mapsto a^2+b^2_1\eta_1+b^2_2\eta_2+c^2\eta_1\eta_2\cr
\endgathered$$
where $a^i,c^i\in R_\odd$ are arbitrary and $b^i_j\in R_\even$
are such that the matrix $(b^i_j)$ is invertible.
The even part of the group is $\GL_2$.
Let us introduce special notations for the elements of the obvious
basis of the Lie superalgebra $\Der({\Lambda[\eta_1,\eta_2]})$ :
\begin{gather}\label{N=2-super-relations}
D_i = {\partial\over\partial\eta_i}, 
\quad D_1^* =\eta_1\eta_2{\partial\over\partial\eta_1},
\quad D_2^* =\eta_2\eta_1{\partial\over\partial\eta_2}
\\
\Theta_i = \eta_i{\partial\over\partial\eta_i}, 
\quad E=\eta_i{\partial\over\partial\eta_2},
\quad F=\eta_2{\partial\over\partial\eta_1}, \quad i=1,2.
\end{gather}
We will call the $D_i$ the {\em differentials}, the $D_i^*$ the {\em homotopies}.
They exhaust the odd basis elements. The even elements $\Theta_i$ will be called
the {\em grading operators}, while $E,F$ will be called the
{\em $\mf{sl}_2$-operators}. Note that
\begin{equation}\label{N=2-homotopies}
[D_1, D_1^*] = \Theta_2, \quad [D_2, D_2^*]=\Theta_1.
\end{equation}

The following classical fact (known to  V. G.  Kac in 1970's, see \cite{Kac-superalgebras}
(3.3.3)),  explains the
special role of the $N=2$ case.

\begin{prop} The group super-scheme 
$G_{\Lambda[\eta_1,\eta_2]}$ is isomorphic to $\SL_{1|2}$.
\end{prop}

\noindent {\sl Proof (sketch):} We have the $1|2$-dimensional superspace
$\Lambda[\eta_1, \eta_2]/\CC \cdot 1\simeq \CC^{1|2}$. The group super-scheme
$G_{\Lambda[\eta_1,\eta_2]}$
acts on this space by linear transformations,
 so we have a morphism of group super-schemes
$G_{\Lambda[\eta_1,\eta_2]}\to\GL_{1|2}$. We then verify directly
that this morphism factors through $SL_{1|2}$. To see that the
resulting morphism $\varphi: G_{\Lambda[\eta_1,\eta_2]}\to\SL_{1|2}$
is an isomorphism, we first verify this on the level
of the underlying even schemes, which are
identified with $GL_2$ for both the source and the target of 
$\varphi$. After this
it remains to verify that $\varphi$ induces an isomorpism
on the level of Lie super-algebras. This is checked directly,
using the above basis in $\Der({\Lambda[\eta_1,\eta_2]})$
and a standard basis in $\mf{sl}_{1|2}$. \qed

\vskip .2cm

Motivated by Proposition \ref{kontsevich-complex-prop}, we
give the following
  
\begin{Defi} An $N=2$ {\em supersymmetric complex} is a super-vector
space $V$ with an action of the group super-scheme 
$G_{\Lambda[\eta_1,\eta_2]}=SL_{1|2}$. 
\end{Defi}

Our goal is now to analyze the structures on an $N=2$ supersymmetric
complex $V$ in more detail. First of all, 
the action on $V$ of the torus $\GG_m\times\GG_m\subset\GL_2$
in the even part of $G_{\Lambda[\eta_1,\eta_2]}$,  gives a bigrading 
$V=\bigoplus_{i,j} V^{ij}$. In other words, the operator $\Theta_1$, (resp. $\Theta_2$)
is equal to $i$ (resp. $j$) on $V^{ij}$. This bigrading is 
compatible with the $\ZZ/2$-grading by parity, i.e.,
 $$V_\even = \bigoplus_{i+j\in 2\ZZ}V^{ij},\quad V_\odd = \bigoplus_{i+j\in 1+2\ZZ}V^{ij},$$
This just expresses the fact that $\Theta_1$ and $\Theta_2$ are even operators.
Next, the operators $D_1$ and $D_2$ define anticommuting differentials
on $V^{\bullet\bullet}$ of square 0 and degrees $(1,0)$ and $(0,1)$ respectively.
This follows from the commutation relations
$$[D_\mu, D_\nu] = 0, \quad [\Theta_\mu, D_\nu] = \delta_{\mu\nu} D_\nu, \quad \mu, \nu=1,2.$$
 which are verified at once from \eqref{N=2-super-relations}. 
 So $V^{\bullet\bullet}$ is, in particular, a double complex. Further, the 
 permutation matrix
 \begin{equation}\label{permutation-matrix}
\begin{pmatrix} 0&1\\1&0\end{pmatrix}\,\,\in\,\,\GL_2(\CC)
 \end{equation}
 identifies $V^{ij}$ with $V^{ji}$ and interchanges $D_1$ and $D_2$. 
 Finally, and most importantly, we have:
 
 \begin{prop}\label{N=2-exactness-proposition}
  Let $V^{\bullet\bullet}$ be an $N=2$ supersymmetric complex.
 Then every row of $V^{\bullet\bullet}$ except, possibly, the $0$th row,
 is exact with respect to $D_1$. Similarly, every column except, possibly,
 the $0$th column, is exact with respect to $D_2$. 
 \end{prop}
 
 \noindent {\sl Proof:} This follows from \eqref{N=2-homotopies}, which means that
 $D^*_\nu$ provides
  a contracting homotopy for $D_\nu$ outside of the $0$th row (for $\nu=1$)
  or the $0$th column (for $\nu=2$). \qed
 
 \vskip .2cm
 
 For a double complex $(C^{\bullet\bullet}, d_1, d_2)$ we denote by $\Tot(C^{\bullet\bullet})$
 its total complex, with differential $d_1+d_2$. 
 
 \begin{cor} Suppose that the bigrading on $V$ is such that $V^{ij}=0$ for
 $i\ll 0$ or $j\ll 0$. Then the complex $\Tot(V^{\bullet\bullet})$
 is quasi-isomorphic to the $0$th row $(V^{\bullet, 0}, D_1)$, as well as to the
 $0$th column $(V^{0\bullet}, D_2)$.
 \end{cor}
 
 \noindent {\sl Proof:} Consider the obvious morphisms of double complexes
 $$V^{\bullet, 0}\buildrel \varphi\over \lla V^{\bullet, \geqslant 0}
 \buildrel\psi\over\lra V^{\bullet\bullet},$$
 with $\varphi$ being surjective and $\psi$ injective. 
 As all the rows of $V^{\bullet,\bullet}$ other than the $0$th row
 are exact, we see that both $\Tot(\Ker(\varphi))$ and $\Tot(\Coker(\psi))$ have
 increasing filtrations with acyclic quotients, whence the statement. \qed
 
 \vskip .3cm
 
 We now consider the particular case where $V^{\bullet\bullet}$ is concentrated in the first
 quadrant, i.e., in the range $i,j\geqslant 0$. Fix $p>0$ and let $V^{pq}_{\cl}\subset V^{pq}$
 be the kernel of $D_1$. We have then the complex
 \begin{equation}
 V^{p,\bullet}_{\cl} \,\,=\,\,\bigl\{ V^{p0}_{\cl}\buildrel D_2\over\hookrightarrow
 V^{p1}_{\cl}\buildrel D_2\over\lra V^{p2}_{\cl}\buildrel D_2\over\lra ...\bigr\}
 \end{equation}
 Note that the very first differential is an injective map as $[D_2, D_2^*]=\Theta_1=p$
 on $V^{pq}$, and we assumed $p>0$. 
 
\begin{prop}\label{V-cl-proposition}
 Let $V^{\bullet\bullet}$ be an $N=2$ supersymmetric
complex concentrated in the first
 quadrant, and $p>0$. Then the complex $ V^{p,\bullet}_{\cl} $ is isomorphic,
 in the derived category, 
 to the cohomological truncation
 $$\bigl(\ten_{\geqslant p+1} (V^{0\bullet}), D_2\bigr)\,\,=\,\,\bigl\{ \Ker(D_2)\hookrightarrow
 V^{0,p}\buildrel D_2\over\lra V^{0, p+1}\buildrel D_2\over\lra V^{0, p+2}
 \buildrel D_2\over\lra...\bigr\}$$
 with the grading normalized so that $\Ker(D_2)$ is in degree $0$.
\end{prop}

\noindent {\sl Proof:} Since the action of the 
permutation matrix \eqref{permutation-matrix}
interchanges the two differentials, it suffices to identify $V^{p,\bullet}_{\cl}$
(up to quasi-isomorphism) with 
$$\bigl(\ten_{\geqslant p+1} (V^{\bullet, 0}), D_1\bigr)\,\,=\,\,\bigl\{ \Ker(D_1)\hookrightarrow
 V^{p,0}\buildrel D_1\over\lra V^{p+1,0}\buildrel D_1\over\lra V^{ p+2,0}
 \buildrel D_1\over\lra...\bigr\},$$
 where $\Ker(D_1)=V^{p,0}_\cl$.
 To achieve this, for each $j\geqslant 0$ consider a similar complex:
 $$\ten_{\geqslant p+1}(V^{\bullet, j})\,\,=\,\,\bigl\{V^{p,j}_{\cl}\lra V^{p,j} \lra V^{p+1,j}\lra ...\bigr\}.$$
 By Proposition \ref{N=2-exactness-proposition}, $\tau_{\geqslant p+1}(V^{\bullet, j})$
 is exact for $j\geqslant 0$. So we consider the double complex 
 \be\label{double-complex-W-equation}
 W^{\bullet\bullet} \,\,=\,\,\bigl\{\ten_{\geqslant p+1}(V^{\bullet, 0})\lra \ten_{\geqslant p+1}(V^{\bullet, 1})\lra
 \ten_{\geqslant p+1}(V^{\bullet, 2})\lra ...\bigr\}
 \ee
 and denote its total complex by $W^\bullet$. Then one edge of $W^{\bullet\bullet}$
 is $\ten_{\geqslant p+1} (V^{\bullet, 0})$, the other edge is $V^{p,\bullet}_{\cl}$
 and all the rows and columns other than these edges are exact. Therefore the projections
 \be\label{two-quasiisomorphisms-equation}
 \ten_{\geqslant p+1} (V^{\bullet, 0})\lla W^\bullet\lra V^{p,\bullet}_{\cl}
 \ee
 of the total complex onto the two edges are quasi-isomorphisms.\qed
 
 \begin{rem}
 Representations of the Lie super-algebra $\Der\,\Lambda[\eta_1, \eta_2]=\mf{sl}(1|2)$
 have attracted a lot of attension. In particular, there is
 a complete classification of finite-dimensional
 irreducible \cite{Bernstein-Leites} and even
 indecomposable \cite{Germony, Leites-indecomposable} representations. 
  In this paper we do not need any more information about these representations
 than what is given by Proposition \ref{V-cl-proposition}. 
 
 \end{rem} 
 
 
 \subsection{The double de Rham complex.}

 Let $X$ be a super-scheme. The super-scheme
 $$\Sc^2X \,\,=\,\,X^{\Lambda[\eta_1, \eta_2]}\,\,=\,\,\ul{\Hom}(\mathbb A^{0|2}, X)$$
 can be seen as the configuration space of an $N=2$ supersymmetric particle
 moving in $X$. The group super-scheme $G_{\Lambda[\eta_1, \eta_2]}$ acts on
 $\Sc^2X$. 
 Denoting by $\varpi^2:\Sc^2X\to X$ and $\varpi: \Sc X\to X$ the projections, we see
 that $\varpi^2_*\Oc_{\Sc^2 X}$ is a sheaf of $N=2$ supersymmetric complexes on $X$. 
 These complexes are concentrated in the first quadrant. 
Viewing $\Sc^2X$ as $\Sc\Sc X$, we can view  $\varpi^2_*\Oc_{\Sc^2 X}$ 
as $\varpi_*\Omega^\bullet_{\Sc X}$, i.e.,
 the de Rham complex of the de Rham complex of $X$.
It has two differentials:
$D_1 = \don_{\on{DR}}^{\Sc X}$, the de Rham differential of $\Sc X$, and $D_2=\on{Lie}_D$,
where $D$ is the vector field on $\Sc X$ corresponding to the de Rham differential
$\don_{\on{DR}}^X$. These differentials are just a part of the structure of an
$N=2$ supersymmetric complex.

\begin{exa} Suppose that $X$ is a purely even smooth
algebraic variety with an etale coordinate system $\phi: X\to \mathbb A^n$, so we have
the  regular functions $x_1, ..., x_n$ on $X$. Then $\Sc X$ has etale coordinates
$x_1, ..., x_n, \xi_1, ..., \xi_n$ with $\xi_i= \don_{\on{DR}}^X (x_i)$
 odd. Accordingly, $\Sc^2 X$
has an  etale coordinate system consisting of $2n$ even coordinates
$x_i$ and $\don_{\on{DR}}^{\Sc X} (\xi_i)$ and $2n$ odd coordinates $\xi_i$ and
 $\don_{\on{DR}}^{\Sc X} (x_i)$, $i=1, ...,n$.
\end{exa}

Let $\Omega^{p,\cl}_{\Sc X}$ be the sheaf of $\don_{\on{DR}}^{\Sc X}$-closed
$p$-forms on $\Sc X$. The direct image $\varpi_* \Omega^{p,\cl}_{\Sc X}$ onto $X$
has the residual grading  and differential
coming from the action of $G_{\Lambda[\eta_2]}$ on $\Sc X$. In local etale
coordinates as above this is the grading assigning  degree 1 to $\xi_i$ and to 
$\don_{\on{DR}}^{\Sc X}(\xi_i)$  and degree 0  to the other generators,
while the differential is $D=D_2$. 
Proposition \ref{V-cl-proposition} implies the following.

\begin{cor}\label{complex-corollary} The complex 
$$\varpi_*(\Omega^{p,\cl}_{\Sc X})^0\buildrel D\over\lra 
\varpi_*(\Omega^{p,\cl}_{\Sc X})^1\buildrel D\over\lra 
\varpi_*(\Omega^{p,\cl}_{\Sc X})^2\buildrel D\over\lra ...$$
 of sheaves on $X$ is quasi-isomorphic to 
$$\Omega^{p, \cl}_X \hookrightarrow
 \Omega^{p}_X \buildrel \don_{\on{DR}}^X\over \lra
 \Omega^{p+1}_X \buildrel \don_{\on{DR}}^X\over \lra
 \Omega^{p+2}_X \buildrel \don_{\on{DR}}^X\over \lra...$$
 where $\Omega^{p,\cl}_X$ is the sheaf of closed differential $p$-forms on $X$. 
\end{cor}


 \section{Super-ind-schemes.}
 
 \subsection{Basic definitions.} We refer to \cite{Artin-Mazur}
 \cite{SGA-IV-I} for general background
 on ind- and pro-objects. 
 By a {\em super-ind-scheme} in this paper we mean an ind-object in $\Schb$
represented as a filtering inductive limit of quasicompact super-schemes and their immersions
\begin{equation}\label{general-form-ind-scheme-eq}
Y\quad =\quad ``\ind_{\a\in A}"\,\, Y^\alpha.
\end{equation}
Alternatively, $Y$  can be identified with the corresponding 
(ind-)representable functor
\begin{equation}
h_Y: \Schb\to \Setb, \quad S\,\,\mapsto\,\,  \ind_{\a\in A}{}^{\Setb}\,\,\, \Hom_{\Schb} (S, Y^\a). 
\end{equation}
We denote by $h^Y$ the covariant functor $\Algb\to\Setb$ given by
$h^Y(R) = h_Y(\Spec(R))$. 

Let $\Ischb$ be the category of super-ind-schemes. It is a general property 
of ind-objects that for a quasi-compact super-scheme $S$
we have $h_Y(S)=\Hom_\Ischb(S,Y)$.

\begin{prop}
Let $X$ be a super-scheme $X$. Consider the object
$$``X" \,\,=\hskip .8cm
\,\,``\hskip -.9cm \ind_{U\subset X \,\,\text{quasicomp.}}\hskip -.9cm " \,\hskip .8cm U
\quad\in\quad \Ischb$$
where $U$ runs over quasi-compact open sub-super-schemes in $X$. 
Associating $X\to ``X"$ defines an embedding of $\SSchb$ into $\Ischb$ as a full
sub-category. \qed
\end{prop}

Our requirement that $X_i$ be quasi-compact follows  \cite{BD2} \S 7.11.
Note that if we defined a super-ind-scheme as simply an ind-object in the category of
all super-schemes, then a non-quasi-compact super-scheme $X$ would be represented by
two ind-objects: $X$ itself and $``X"$, which are not isomorphic in general.

\vskip .2cm

Let $f: X\to Y$ be a morphism of super-ind-schemes. We say that $f$ is {\em formally smooth},
if the induced 
morphism $h^X\to h^Y$ of contravariant functors  $\Affb\to\Setb$ is formally smooth in
the sense of Definition \ref{formal-smoothness-definition}(a). 
The even  and reduced parts of an super-ind-scheme $Y$ as
 in \eqref{general-form-ind-scheme-eq}
 is defined by 
\be
Y_{\on{even}}\,\,=\,\, ``\ind_{\a\in A}" \,\,\,Y^\a_{\on{even}},\quad
Y_{\red}\,\,=\,\, ``\ind_{\a\in A}" \,\,\,Y^\a_{red}. 
\ee
Similarly, let $\oen$ be a finite dimensional local $\CC$-super-algebra.
Using Proposition \ref{near-points-proposition},
 we extends the functor $X\to X^\oen$ to ind-schemes
by
 \begin{equation} Y^\oen\,\,
 =\,\,
 ``\ind_{\a\in A}"\,\,(Y^\a)^\oen.
 \end{equation}

 \begin{exa} Let $B$ be a super-scheme, $I$ a finite set, and $B^I$ the $I$th Cartesian
 power of $B$. A morphism $u: S\to B^I$ is thus the same an as $I$-tuple of morphisms
 $u_i: S\to B$, $i\in I$. Denoting $\Delta\subset B^I$ the small diagonal $\{(b,b,...,b)\}$
 and by $\Ic_\Delta\subset\Oc_{B^I}$ its sheaf of ideals, we can view the 
 formal neighbourhood of $\Delta$ in $B^I$ as an ind-scheme
 $$B^{[I]}\,\,=\,\,``\ind_{n\geqslant 0}" \, B^{[I]}_n, \quad 
  B^{[I]}_n \,\,=\,\,
 \Spec_{B^I}\, \bigl(\Oc_{B^I}/\Ic_{\Delta}^{n+1}\bigr).$$ 
 A morphism from a super-scheme $S$ into $B^{[I]}$ is the same as an $I$-tuple of
 morphisms $u_i: S\to B$ as above but with the condition that any two $u_i, u_j$
 are infinitesimally near, in the sense of Definition \ref{inf-near-morphisms-def}.
 Note that for a 1-element set $I$ we have $B^I=B^{[I]}=B$. 
 
Further, any map $p: J\to I$ of finite sets induces a morphism of schemes
$p^*: B^I\to B^J$ and a morphism of ind-schemes $[p]^*: B^{[I]}\to B^{[J]}$.
If $p$ is injective, then $p^*$ and $[p]^*$ are coordinate projections;
if $p$ is surjective, then  $p^*$ and $[p]^*$ are diagonal embeddings. 
 
 \end{exa}
 
 We now discuss the concept of an integrable connection, following the approach of
 Grothendieck \cite{Grothendieck-crystals}, see also \cite{BD1} (3.4.7). 
 
 Let $B$ be a super-scheme, and $E\to B$ be a super-ind-scheme over $B$. 
 For a morphism of super-schemes $u: S\to B$ we denote by $u^*E= E\times_BS\to S$
 the pullback of $E$. 
 
\begin{prop}\label{int-connections-conditions-prop}
 For a given $E\to B$ as above, the following systems of  data (1) and (2) 
are in a bijection:

(1) For each super-scheme $S$ and each pair of infinitesimally near morphisms
$u,u': S\to B$, an isomorphism $M_{S,u,u'}: u^*E\to u'{}^*E$ of super-ind-schemes over $S$,
satisfying the following conditions:

(1a) Transitivity: for each three infinitesimally near morphisms
$u,u', u'': S\to B$, we have
$$M_{S,u,u''}\,\,=\,\,M_{S,u',u''}\circ M_{S,u,u'}.$$

(1b) Compatibility with restrictions: for any $u,u'$ as above and any morphism
$v: S'\to S$, we have 
$$M_{S', uv, u'v}\,\,=\,\,v^*M_{S,u,u'}.$$

(2) For each nonempty finite set $I$, an super-ind-scheme $E_I\to B^{[I]}$ such that 
$E_I=E$ for any 1-element $I$, and for any map $p: J\to I$ we have an isomorphism
$\alpha_p: E_I\to [p]^*E_J$, these isomorphisms compatible with compositions of maps. 
\end{prop}  

We will call a datum of either type an {\em integrable connection}
on $E$ along $B$.

\vskip .2cm
 
 \noindent {\sl Proof:} Given a datum of type (2), any infinitesimally near
 $u,u':S\to B$ give a morphism $(u,u'): S\to B^{[\{1,2\}]}$. On the other hand, 
 $B^{[\{1,2\}]}$ is the formal neighborhood of the diagonal in $B\times B$, and
 the isomorphisms $\alpha_{i_1}, \alpha_{i_2}$ corresponding to the 
 maps $i_1: \{1\}\hookrightarrow \{1,2\}$, $i_2: \{2\}\hookrightarrow \{1,2\}$
 identify
 $E_{[\{1,2\}]}$ with the pullback of $E$ via the two projections 
 $[i_1]^*, [i_2]^*: B^{[\{1,2\}]}\to B$, whence the isomorphism $M_{S,u,u'}$.
 Transitivity follows from considering the morphism  $(u,u',u''): S\to B^{[\{1,2,3\}]}$.
 Compatibility with restrictions follows because the morphism
 $(uv,u'v): S'\to B^{[\{1,2\}]}$ is the composition of $(u,u')$ and $v$.
 
 \vskip .2cm
 
 Conversely, let a datum of type (1) be given. To construct $E_I$, we fix $n\geqslant 0$
 and take $S=B^{[I]}_n$. The coordinate projections $p_i: S\to B$, $i\in I$,
 are infinitesimally close to each other so the ind-schemes $p_i^*E\to B^{[I]}_n$
 are canonically identified with each other via tha $M$-isomorphisms.
  We can say that we have one ind-scheme
 $E_{I,n}$, identified with them all. When $n$ increases, these $E_{I,n}$ form an
 filtering inductive system of super-ind-schemes and closed embeddings, so
 their limit $E_I$ is a well defined object of $\Ischb$. 
  The remaining verifications are left to the reader. 
 \qed

 
 \subsection{Functions and forms on
 super-ind-schemes.}\label{functions-forms-superschemes-subsection}
 As in  \cite{Haboush} and  \cite{KV2}, Sect.~2, to any super-ind-scheme $Y$ as in
 \eqref{general-form-ind-scheme-eq},
 we associate a topological space 
\be \ul{Y}\,\,=\,\,
\ind_{\a\in A}{}^\Topb \,\,\, \ul{Y}^\a.
\ee
Let $i_\a: \ul{Y}^\a\to \ul{Y}$ be  the canonical embedding.
We then have a sheaf of super-commutative pro-algebras $\Oc_Y$ on $\ul{Y}$
\be 
\Oc_Y\,\,=
\,\, \pro_{\a\in A}\,\,(i_\a)_*\Oc_{Y^\a},
\ee
 which we can consider as a sheaf of topological algebras.
We define the sheaves of differential forms in a similar way: 
\be 
\Omega_Y^p\,\,=\,\, \pro_{\a\in A}\,\,(i_\a)_*\,\Omega^p_{Y^\a}.
\ee
Note that
\be \Omega^\bullet_Y\,\,=\,\, \varpi_*\Oc_{\Sc Y},
\ee
where $\pi:\Sc Y\to Y$ is the natural projection.

 
 \section{The formal loop space of a super-manifold.}
 
 \subsection{ Nil-Laurent series}
 
For a super-commutative ring $R$ we denote by 
$R\lb t \rb^\sqr$ the subring of $R\lb t\rb$ consisting of  Laurent series 
$\sum_{i\gg -\infty}^\infty a_i t^i$
such that $a_i$ is nilpotent for $i<0$. We proved in \cite{KV1}, Prop. 1.3.1, that
if $R$ is a commutative local ring, then so is $R\lb t\rb^\sqr$. 
We need the following version of this. 

 \begin{lem}\label{nil-local-ring-lemma}
 Let $S$ be a super-scheme. Then $\Oc_S\lsb t\rsb$ and
$\Oc_S\lb t\rb^\sqr$ are sheaves of super-commutative local rings.
 \end{lem}

\noindent{\sl Proof :}
It is enough to assume that $S=\Spec(R)$. Let $\pen\in\underline{\Spec}(R)$,
i.e., $\pen\subset R_0$ is a prime ideal.
We first treat the case of $\Oc_S\lsb t\rsb$. The stalk of this sheaf at $\pen$
is the ring 
$$\Oc_S\lsb t\rsb_\pen\,\,\,=\,\,\,\ind_{U\ni\pen}\,\Oc(U)\lsb t\rsb
\,\,\,=\,\,\,
R\lsb t\rsb\bigl[(R-\pen)^{-1}\bigr].$$ 
We claim that 
$$\pen'\,\,\,=\,\,\, \biggl\{ b^{-1} \sum_{n=0}^\infty a_nt^n;\quad a_n\in R,\,\,\,a_0\in\pen,
\,\,\,b\in R-\pen
\biggr\}$$
is the maximal ideal in
$R\lsb t\rsb\bigl[(R-\pen)^{-1}\bigr]$, i.e., any element not in $\pen'$ is invertible.
This is obvious by using the geometric series and inverting $a_0\notin\pen$. 
 
\vskip .2cm

Consider now the case of $\Oc_S\lb t\rb^\sqr$. As before, we have
$$\Oc_S\lb t\rb^{\sqr}_\pen \,\,\, =\,\,\, R\lb t\rb^\sqr\bigl[(R-\pen)^{-1}\bigr].$$ 
We define
\begin{equation}\label{tilde-p-equation}
\tilde\pen\,\,\,=\,\,\,\biggl\{b^{-1}\sum_{n\gg-\infty}^\infty a_nt^n;\quad a_n\in R,
\,\,\,
a_{<0}\in\sqrt{R},\,\,\,a_0\in\pen,\,\,\,b\in R-\pen\biggr\}
\end{equation}
and claim that it is the maximal ideal in
$\Oc_S\lb t\rb^\sqr_\pen.$ 
Indeed, the fact that $\tilde\pen$ is an ideal is obvious.
On the other hand, if 
$u(t)\in\Oc_S\lb t\rb^\sqr_\pen-\tilde\pen$, 
then we write $u(t)$ as the sum
$$u(t)=u_-(t)+a_0b^{-1}+u_+(t)$$ where $u_\pm(t)$ is the sum
of the terms with $\pm n>0$. Now, $u_-(t)$ is nilpotent, 
$a_0b^{-1}$ is invertible in $R\bigl[(R-\pen)^{-1}\bigr]$, and $u_+(t)$ is topologically
nilpotent. So the invertibility follows in the same way as in 
\cite{KV1}, Prop.~1.3.1.
\qed 

\vskip .3cm

As in \cite{KV1} (1.6), denote by $\Eb$ the set of sequences
\be\label{epsilon-set-equation}
\epsilon = (\epsilon_{-1}, \epsilon_{-2}, ...), \quad \epsilon_j\in \ZZ_+, \quad
\epsilon_j=0, \quad j\ll 0.
\ee
It is equipped with a natural partial order such that $\epsilon \leq\epsilon'$
if $\epsilon_j\leq\epsilon'_j$ for all $j$. For a super-commutative algebra $R$
we define the subset
\be 
R\lb t\rb^\sqr_\epsilon \,\,\, = \,\,\, \biggl\{ \sum_{n\in \ZZ} a_n t^n \biggl| \,\,
a_n^{1+\epsilon_n} =0, \,\, n<0\biggr\}.
\ee
Thus series from this set have both the number of negative coefficients
and their order of nilpotency bounded.

\begin{prop}\label{finite-subalgebra-R((t))-proposition}
Any finitely generated subalgebra $A$ in $R\lb t\rb^\sqr$
is contained in $R\lb t\rb^\sqr_\epsilon$ for some $\epsilon$. 
\end{prop}

\noindent {\sl Proof:} Let $f_1, ..., f_r$ be generators of $A$,
which we can assume to be homogeneous with respect to the
$\ZZ/2$-grading. Write $f_i = f_{i,+}+f_{i,-}$, where $f_{i,+}\in R\lsb t\rsb$,
while $f_{i,-}$ is the sum of the terms with negative powers of $t$.
Then, each $f_{i,-}$ is nilpotent. This implies that among the infinite
number of monomials
$$
f_-^m:= f_{1,-}^{m_1}f_{2,-}^{m_2}\, ...\, f_{r,-}^{m_r},\quad m=(m_1, ..., m_r),
 \quad m_i \geqslant 0,
$$
only finitely many are non-zero. Let $m^{(1)}, ..., m^{(s)}$ be the 
exponents of all the nonzero ones. 
Look now at similar monomials
$f^n=f_1^{n_1}...f_r^{n_r}$ formed out of the $f_i$. They form a spanning set for $A$.
On the other hand, expanding them using $f_{i}=f_{i,+}+f_{i,-}$ and the
binomial formula, we find that each $f^n$ can be expressed as
$$ f^n\,\,=\,\,\sum_{\nu=1}^s F_\nu^n f_-^{m^{(\nu)}},\quad F_\nu^n\in R\lsb t\rsb.$$
The finitely many monomials $f_-^{m^{(\nu)}}\in R\lb t\rb^\sqr$
clearly admit $N,d\geqslant 0$ with the following properties.
 First,  all the $f_-^{m^{(\nu)}}$
 have zero coefficients at $t^j$, $j<-N$. 
Second, all the coefficients of these $f_-^{m^{(\nu)}}$
at monomials with $t^{j}$, $-N\leqslant j\leqslant -1$, are nilpotent of degree $d+1$.
Look now at elements of the form
  $F f_-^m$ with $F\in R\lsb t\rsb$. Each of them clearly satisfies the first
  property: the order of pole is still bounded by $N$. As for the second
  property, each coefficient of $Ff_-^m$ at each negative power of $t$
  is a sum of at most $N-1$ summands, each nilpotent of degree $d+1$.
  This implies that there is $d'$ depending only on $d$ and $N$
   such that each coefficients of each
   $Ff_-^m$ at each negative power of $t$, is nilpotent of
  degree $d'+1$. 
  This means that  $A\subset R\lb t\rb^\sqr_\epsilon$, where
$\epsilon$ is such that $\epsilon_{-1}=...=\epsilon_{-N}=d'$, and $\epsilon_i=0$
for $i<N$.\qed


\subsection{Basics on $\Lc^0X$ and $\Lc X$.}
Let $X$ be a super-scheme.
We define the super-scheme
$$\Lc^0_nX=X^{\CC[t]/t^{n+1}}.$$
For different $n$ the $\Lc^0_nX$ form a projective system of affine morphisms
of super-schemes. We define the super-scheme $\Lc^0X$ to be the projective 
limit of this system, and call it the {\em super-scheme of formal arcs} in $X$.
Compare with \cite{DL}.

 \begin{prop}\label{formal-arcs-functor-prop}
For any super-commutative ring $R$ and, more generally,
for any super-scheme $S$ we have
$$\gathered
\Hom_{\SSchb} ({\Spec}(R), \Lc^0 X) \quad = \quad \Hom_{\SSchb}(
{\Spec}(R\lsb t\rsb ), X), \cr
\Hom_{\SSchb}(S, \Lc^0 X) \quad =\quad \Hom_{\Sspb}\bigl((\ul{S}, \Oc_S\lsb t\rsb ), \, X
\bigr).\endgathered$$
\end{prop}

This was asserted for schemes
in \cite{KV1}, Prop. 1.2.1(b) but with an incorrect proof
(the first equality in Lemma 1.2.3 of {\em loc. cit.} does not
hold in general). Here we supply the proof. 

 \vskip .2cm

\noindent {\sl Proof of \eqref{formal-arcs-functor-prop}:} 
Note that
if $S$ is any super-scheme, then, by Proposition 
\ref{near-points-proposition} applied
to $\oen=\CC[t]/t^{n+1}$ we have
$$\Hom_{\SSchb}(S, \Lc^0_n X) \quad =\quad \Hom_{\Sspb}\bigl((\ul{S}, \Oc_S[t]/t^{n+1}), X
\bigr).$$
Next, we have
$$\Oc_S\lsb t\rsb  \quad = \quad \pro_{n\geqslant 0}\,\, \Oc_S[t]/t^{n+1}$$
in the category of sheaves of local rings on $\ul{S}$, so
$$(\ul{S}, \Oc_S\lsb t\rsb ) \quad =\quad \pro_{n\geqslant 0}{}^{\Sspb} \, \bigl(\ul{S}, \Oc_S[t]/t^{n+1}
\bigr),$$
and therefore
$$\Hom_{\Sspb}\bigl((\ul{S}, \Oc_S\lsb t\rsb ), X\bigr) \quad = \quad \pro_{n\geqslant 0}
\, \Hom_{\Sspb}\bigl((\ul{S}, \Oc_S[t]/t^{n+1}), X\bigr).$$
Note that $\SSchb$ is a full subcategory in $\Sspb$, 
so $\Hom$ on the right hand side
can be taken in either category. Now the fact that
$$\Lc^0 X \quad =\quad \pro_{n\geqslant 0} {}^{\SSchb} \, \,\, \Lc^0_n X$$
implies that 
$$\Hom_{\SSchb}(S, \Lc^0 X)\,\,\, = \,\,\,
 \pro_{n\geqslant 0} \,\, \Hom_{\SSchb}(S, \Lc^0_n X)\,\,\, = \,\,\,
\Hom_{\Sspb}\bigl( (\ul{S}, \Oc_S\lsb t\rsb ), X\bigr),$$
as claimed. \qed

\vskip .3cm

As in [KV1] we define the functor $\lambda_X:\SSchb\to\Setb$ as follows :
\be 
\lambda_X(S)\,\,\,=\,\,\, \Hom_\Sspb\bigl((\ul{S},\Oc_S\lb t\rb^\sqr),X\bigr).
\ee
\begin{prop}\label{formal-loop-functor-prop}
(a) If $X=\Spec(A)$, $S=\Spec(R)$ are affine super-schemes, then
$$\lambda_X(S)=\Hom_\Algb\bigl(A, R\lb t\rb^\sqr\bigr).$$

(b) For any super-scheme
$X$ of finite type the functor $\lambda_X$ is representable by a 
super-ind-scheme $\Lc X$, and 
$\Lc X=\ind{}_{U\subset X\,\text{affine}}\Lc U$ in the 
category of super-ind-schemes. 
\end{prop}

\begin{rem}

In \cite{KV1}, Prop.~1.4.5, we claimed (with an incorrect proof, based on
 erroneous Lemma 1.4.3(a)), that the analog of \eqref{formal-loop-functor-prop}(a)
  holds for any
 $X$ of finite type. In fact, this stronger statement is unnecessary, and 
 \eqref{formal-loop-functor-prop}(a)
 is sufficient to establish \eqref{formal-loop-functor-prop}(b) and all the properties of $\Lc X$ claimed in \cite{KV1}.
 
 \end{rem}

\vskip3mm

\noindent{\sl Proof of Proposition \ref{formal-loop-functor-prop}(a):}
 Let $f\in\lambda_X(S)$, i.e.,
$$f=(f_\flat,f^\sharp):\,\, \bigl(\underline{\Spec}(R),\Oc_{\Spec(R)}\lb t\rb^\sqr\bigr)
\,\,\,\lra\,\,\,
\bigl(\underline{\Spec}(A),\Oc_{\Spec(A)}\bigl)$$
is a morphism
of super-spaces. Thus
$f_\flat:\underline{\Spec}(R)\to\underline{\Spec}(A)$
is a morphism of topological spaces, and
$$f^\sharp:\,\, f_\flat^{-1}\Oc_{\Spec(A)}
\,\,\,\lra\,\,\,
\Oc_{\Spec(R)}\lb t\rb^\sqr$$
is a morphism of sheaves of super-commutative local rings.
It induces a morphism of rings 
$$\varphi=\Gamma(f^\sharp):\,\, A=\Gamma(\underline{\Spec}(A),\Oc_{\Spec(A)})
\,\,\lra\,\,
R\lb t\rb^\sqr=\Gamma\bigl(\underline{\Spec}(R),\Oc_{\Spec(R)}\lb t\rb^\sqr\bigr)$$
and so a morphism of super-schemes
$$g:\,\,\Spec R\lb t\rb^\sqr \,\,\lra\,\, \Spec(A).$$
So it is enough to prove:

\begin{lem}
The correspondence $f\mapsto\varphi$ gives a bijection
$$\Phi:\,\, \lambda_X(S)\,\,\,\to\,\,\,\Hom_\Algb\bigl(A,R\lb t\rb^\sqr\bigr).$$
\end{lem}

\noindent{\sl Proof :}
We construct the inverse map
$$\Psi:\,\, \Hom_\Algb\bigl(A,R\lb t\rb^\sqr\bigr)\,\,\,\to\,\,\,\lambda_X(S).$$
We have a morphism of super-spaces
$$h=(h_\flat,h^\sharp):\,\,
\bigl(\underline{\Spec}(R),\Oc_{\Spec(R)}\lb t\rb^\sqr\bigr)\,\,\,\lra\,\,\,
\Spec (R\lb t\rb^\sqr)$$
with the map of topological spaces $h_\flat$ defined as the composition
$$\underline{\Spec}(R)\,\,\,=\,\,\,\underline{\Spec}(\overline {R})
\,\,\,{\buildrel u\over\lra}\,\,\,
\underline{\Spec}(\overline{ R}\lsb t\rsb )
\,\,\,{\buildrel v\over\lra}\,\,\,
\underline{\Spec}(R\lb t\rb^\sqr).
$$
Here we have denoted $\overline{ R}=R/\sqrt{R}$, and $u$ is induced by the evaluation
homomorphism
$$\overline{ R}\lsb t\rsb\,\,\,\lra\,\,\, \overline{ R} \,\,\,=\,\,\,
\overline{ R}\lsb t\rsb/t\overline{ R}\lsb t\rsb,$$
while $v$ is induced by the termwise factorization by $\sqrt{R}$: 
$$R\lb t\rb^\sqr\,\,\,\lra\,\,\,\overline{ R}\lb t\rb^\sqr \,\,\,=\,\,\,
\overline{ R}\lsb t\rsb.$$
The morphism $h^\sharp$ is induced by the inclusions
$$R\lb t\rb^\sqr [1/b]\,\,\,\subset\,\,\,(R[1/b])\lb t\rb^\sqr,\quad b\in R.$$
Given $\varphi:A\to R\lb t\rb^\sqr$, it induces a morphism of super-schemes
$$g:\Spec R\lb t\rb^\sqr \,\,\,\lra \,\,\, \Spec(A),$$
and we define $f=\Psi(\varphi)$ to
be the composition
$$\Psi(\varphi)\,\,\, = \,\,\, gh:\,\, 
\bigl(
\underline{\Spec}(R), \Oc_{\Spec(R)}\lb t\rb^\sqr \bigr) \,\,\,\to \,\,\,
\Spec(A) = X.$$

We now claim that the maps $\Phi$ and $\Psi$ are inverse to each other. Indeed,
the equality $\Phi\Psi=\Id$ is obvious, it follows from the fact that 
$\Gamma(h^\sharp)$ is the identity of $R\lb t\rb^\sqr$.
 
Let us prove that $\Psi\Phi=\Id$.
The proof is analogous to the classical proof that a morphism
of affine schemes is the same as a 
homomorphism of the corresponding rings.
So let $f=(f_\flat,f^\sharp)\in\lambda_X(S)$, and
$g=(g_\flat,g^\sharp)=\Phi(f)$.
By construction
$$\Psi(g)\,\,\, =\,\,\, \bigl(g_\flat h_\flat,g_\flat^{-1}(h^\sharp)g^\sharp\bigr).$$
Let us prove the equality of maps $g_\flat h_\flat=f_\flat$,
leaving the other equality to the reader.

Let $\pen\in\underline{\Spec}(R)$, so $\pen\subset R_0$ is a prime ideal.
By definition of $g_\flat$ the equality $f_\flat(\pen)=g_\flat h_\flat(\pen)$ is
equivalent to
\be\label{LX-lemma-equation}
f_\flat(\pen)=\varphi^{-1}h_\flat(\pen),
\ee
where $\varphi=\Gamma(f^\sharp)$ fits into the commutative diagram
 $$
\begin{matrix}
A&{\buildrel\varphi\over\lra}&R\lb t\rb^\sqr\cr
\big\downarrow&&\big\downarrow\cr
A_{f_\flat(\pen)}&{\buildrel f_\pen^\sharp\over\lra}&
(\Oc_{\Spec R}\lb t\rb^\sqr)_\pen.
\end{matrix}
$$ 
The vertical maps in this diagram are obtained by taking the stalks, and the map
$f^\sharp_\pen$ is a local homomorphism of local rings.
We use the notation \eqref{tilde-p-equation} for the maximal ideal $\tilde\pen$
in $(\Oc_{\Spec R}\lb t\rb^\sqr)_\pen$.
Let $\tilde{\tilde\pen}$ be its inverse image in $R\lb t\rb^\sqr$.
Explicitly, we have
$$\tilde{\tilde\pen}\,\,\,=\,\,\, \biggl\{\sum_{n\gg -\infty}^{+\infty}a_nt^n\in R\lb t\rb^\sqr;
\quad 
a_0\in\pen\biggr\}.$$
Since the diagram above commutes and $f_\pen^\sharp$ is a local homomorphism, 
we have
$\varphi^{-1}(\tilde{\tilde\pen})=f_\flat(\pen)$.
So to prove \eqref{LX-lemma-equation} we need to show that $h_\flat(\pen)=\tilde{\tilde\pen},$
which is obvious.

This ends the proof of Proposition 
\ref{formal-loop-functor-prop}, part (a).
The proof of part (b) is then achieved as in \cite{KV1}.
Indeed, for $X$ affine, part (a) implies that $\lambda_X$ is represented by the 
formal neighbhorhood of $\Lc^0X$ in $\tilde\Lc X$, see \cite{KV1}, p.~219.
For general $X$ of finite type, $\Lc X$ is glued from $\Lc U$, $U\subset X$ 
affine, as in \cite{KV1}, Prop.~1.4.6.\qed

\vskip .3cm

Recall the De Rham spectrum functor $\Sc$ from Subsection \ref{N=1-supersymmetry-subection}.

\begin{prop}\label{LSX=SLX-proposition}
For any super-scheme $X$ of finite type we have an isomorphism
of super-ind-schemes $\Lc \Sc X=\Sc\Lc X$.
\end{prop}

\noindent{\sl Proof :}
Both super-ind-schemes represent the same functor
$$S\,\,\mapsto\,\,\Hom_\Sspb\bigl((\ul{S},\Oc_S\lb t\rb^\sqr[\eta]),(\ul{X},\Oc_X)\bigr),$$
where $\eta$ is an odd generator, so that $\eta^2=0$.
Indeed, for any super-commutative ring $R$ we have
$$(R[\eta])\lb t\rb^\sqr\,\,=\,\, R\lb t\rb^\sqr[\eta].$$

\qed


\subsection{$\Lc X$ and loco-modules of Borisov.}

By construction that there are morphisms
\be 
X\buildrel\pi\over\longleftarrow \Lc^0 X\buildrel i\over\longrightarrow
\Lc X,
\ee
where $\pi$ is affine and $i$
realizes $\Lc X$ as a formal thickening of $\Lc^0 X$. They are induced by the
obvious
morphisms of sheaves of local rings on any super-scheme $S$:
$$\Oc_S \,\,\,\lla\,\,\, \Oc_S[[t]]
\,\,\,\hookrightarrow \,\,\, \Oc_S((t))^\sqr.$$
We are going to describe explicitly $\pi_*\Oc_{\Lc^0 X}$,
which is a quasicoherent sheaf of $\Oc_X$-algebras,
and $\pi_*\Oc_{\Lc X}$, which is a sheaf of pro-$\Oc_X$-algebras.

\vskip .2cm

Let $A$ be a super-commutative
algebra. Specializing \eqref{A^o-generators-eq} to the particular
case of $\oen = \CC[t]/t^{n+1}$ and of the basis of $\oen$ formed by $1,t, ..., t^n$,
we find:

\begin{cor}\label{borisov-finite-corollary} The super-scheme $\Lc_n^0(\Spec A)$ is identified with
$\Spec(A^{\CC[t]/t^{n+1}})$, where $A^{\CC[t]/t^{n+1}}$ the super-commutative algebra
generated by the symbols $a[m]$, $0\leqslant m\leqslant n$ such that the $\ZZ/2$-degree of $a[m]$
is the same as that of $a$, and which are subject to the relations:
$$(a+b)[m] \,\,= \,\, a[m] + b[m], \quad (\lambda a) [n]
\,\,= \,\, \lambda(a[n]), \,\,\,\lambda\in\CC;
\leqno \eqref{borisov-finite-corollary}(a) $$
$$1[m] =0, \quad m\neq 0; \leqno \eqref{borisov-finite-corollary}(b)$$
$$(ab)[m] \,\,\, = \,\,\, \sum_{i+j=m} a[i] \cdot b[j]. \leqno 
\eqref{borisov-finite-corollary}(c)$$
\qed
\end{cor}

We denote 
 \be\label{A[[t]]-projection-eq}
 A^{\lsb t\rsb}\,\,=\,\,\ind_{n\geqslant 0}\,\, A^{\CC[t]/t^{n+1}}.
 \ee
This algebra can be defined by generators $a[m]$ given for all $m\geqslant 0$ subject to
the same relations as in  \eqref{borisov-finite-corollary}(a)-(c).
 Note that we have an embedding of algebras
\be
A \,\hookrightarrow \, A^{[[t]]}, \quad a\mapsto a[0]. 
\ee

By applying  the limit construction  (inductive for algebras, projective for schemes) 
to  the above corollary and to Proposition \ref{near-points-proposition}, we obtain:

\begin{prop}
(a) If $X=\Spec(A)$, then $\Lc^0 X = \Spec(A^{[[t]]})$, with the projection
$\pi$ induced by \eqref{A[[t]]-projection-eq}.

(b) If $S\subset A$ is a multiplicative subset, then $(A[S^{-1}])^{[[t]]} =
A^{[[t]]}[S^{-1}]$. In particular, for any super-scheme $S$ the sheaf
$\Oc_X^{[[t]]}$ is quasicoherent.

(c) We have an identification
$$\pi_*\Oc_{\Lc^0 X} \,\,\, = \,\,\,\Oc_X^{[[t]]}. $$

\end{prop}

To be precise, (a) follows from Corollary
\ref{borisov-finite-corollary} since projective limits of affine super-schemes correspond to
inductive limits of algebras. Part (b) follows from Proposition 
\ref{near-points-proposition}(c) since localization commutes with inductive
limits. Finally, part (c) follows from part (a).\qed

 \vskip .3cm
 
For each sequence $\epsilon\in\Eb$ as in \eqref{epsilon-set-equation}
let $A^{\lb t\rb }_\epsilon$ be the
algebra with generators $a[n]$ for $a\in A$ and $n\in \ZZ$ (arbitrary integers),
subject to the relations
\be\label{A((t))-epsilon-equation}
a[n]^{1+\epsilon_n} \, = \, 0, \quad a\in A, \, n<0,
\ee
together with the relations 
identical to \eqref{borisov-finite-corollary}(a)-(c) but with $n, i,j\in \ZZ$.
Note that \eqref{A((t))-epsilon-equation} implies that $a[n]=0$ for any $a$ and $n\ll 0$, 
so the sum
in \eqref{borisov-finite-corollary}(c) remains finite.

For $\epsilon\leq\epsilon'$ 
we have a surjection of algebras $A^{\lb t\rb}_{\epsilon'}\to
A^{\lb t\rb}_\epsilon$, and we define the pro-algebra
\be A^{\lb t\rb} \,\,\, = \,\,\, \pro_{\epsilon\in\Eb}\,\, A^{\lb t\rb}_\epsilon. 
\ee
Recall, see Subection \ref{functions-forms-superschemes-subsection},
 that every super-ind-scheme $Y$ gives a topological space
$\ul{Y}$ and a sheaf $\Oc_Y$ over $\ul{Y}$ of pro-supercommutative rings.
In particular, if $Y=\Lc X$ then $\ul{Y}=\ul{\Lc^0 X}$, so that we have a sheaf
$\Oc_{\Lc X}$ over $\ul{\Lc^0X}$.

\begin{prop}\label{borisov-description-LX-prop}
(a) Let $X=\Spec(A)$ be an affine super-scheme of finite type
(i.e., $A$ is finitely generated as an algebra). Then
$$\Lc X \,\,\, = \,\,\,\operatorname{Spf}\,A^{\lb t\rb} \,\,\, := \,\,\,
``\ind_{\epsilon\in \Eb}" \, \Spec \, A^{\lb t\rb}_\epsilon.$$

(b) If $X$ is any super-scheme of finite type, then
we have an identification of sheaves of pro-algebras on $X$
$$\pi_*\Oc_{\Lc X} \, = \, \Oc_X^{\lb t\rb}.$$
\end{prop}

\noindent {\sl Proof:} (a) By definition, for any supercommutative algebra $R$ we have
the first of the following two equalities:
$$\Hom(\Spec\, R, \Lc X)\,\,=\,\,\Hom_{\Algb} (A, R\lb t\rb^\sqr)\,\, =\,\,
\ind_{\epsilon\in\Eb}\,\, \Hom_{\Algb} (A, R\lb t\rb^\sqr_\epsilon)
.$$
The second equality is a consequence of Proposition \ref{finite-subalgebra-R((t))-proposition},
since $A$ is assumed finitely generated. It remains to notice that 
$$\Hom_{\Algb} (A, R\lb t\rb^\sqr_\epsilon)
\,\, =\,\,
\Hom (\Spec\, R, \Spec\,  A^{\lb t\rb}_\epsilon).$$
This proves (a) since a super-ind-scheme is uniquely determined by the functor
 it represents on affine super-schemes. Part (b), being a local statement,  follows from (a) \qed
 
 \begin{rem} The above considerations are very similar to the work of Borisov 
 \cite{Borisov}.
 In particular, his ``loco-modules" can be understood as sheaves of discrete
 modules over the sheaf of  topological (or pro-) algebras $\Oc_X^{\lb t\rb}$,
 i.e., as certain sheaves on the ind-scheme $\Lc X$. 
 
 \end{rem}
 
\begin{exa}
Let $X = \mathbb A^N$, so $A \,= \,
\CC[a_1, ..., a_N]$. For $i=1, ..., N$ and $n\in \ZZ$ let
$b^i_n = a_i[n]\in A^{\lb t\rb}$. Thus the $b^i_n$ are the components
of $N$ indeterminate power series
$$a_i(t) \,\,\, = \,\,\,\sum_{n\gg -\infty} b^i_n t^n$$
forming a point of $\Lc\mathbb A^N$. We have then:
$$A^{\lsb t\rsb} \,\,\, = \,\,\,\CC\bigl[ b^i_n,\,\,\, i=1, ..., N, \,\, n\geqslant 0\bigr];$$
$$A^{\lb t\rb} \,\,\, = \,\,\, \pro_{m>0} \,\,
\CC\big[ b^i_n,\,\,\, i=1, ..., N, \,\, n\geqslant 0\bigr] \bigl[\bigl[ b^i_n,\,\,\,
i=1, ..., N,\, n=-m, ..., -1\bigr]\bigr].$$
The case when $X=\mathbb A^{d_1|d_2}$ is a super-affine space, is considered
similarly: we have even and odd coordinates $a_1, ..., a_N$, $N=d_1+d_2$,
and use the convention of Example \ref{affine-space-examples}(b) for super-polynomial
rings. 

\end{exa}

\begin{rem}\label{L-n-epsilon-remark}
Assume now that $X$ is a smooth algebraic super-variety of dimension
$d_1|d_2$.
If $U\subset X$ is a Zariski open
set admitting an etale map $\phi: U\to \mathbb A^{d_1|d_2}$, then
$\Lc U \subset \Lc X$ is open.
Then $\Lc U$ admits a representation as the limit of a Cartesian
ind-pro-system as in \cite{KV1} :
$$\Lc U \,\,\, = \,\,\, ``\ind_{\epsilon\in\Eb}" \,\, \pro_{n\geqslant 0}\,\, \Lc^\epsilon_n(\phi).
$$
Here  
the scheme $\Lc^\epsilon_n(\phi)$ is defined as follows.
First, consider the case when $U=\mathbb A^{d_1|d_2}$ 
with (even and odd) 
coordinates $a_1, ..., a_N$, $N=d_1+d_2$, and $\phi = \Id$. In this case
$$\Lc^\epsilon_n(\Id) \,\,\, = \,\,\, \Spec \,\, 
\CC\biggl[a_i[l];  -N\leqslant n\leqslant l\biggr]\biggr/ \biggl(
(a_i[l])^{1+\epsilon_l}; l<0\biggr), $$
where $N$ is any number such that $\epsilon_l=0$ for $l<-N$.
This is a super-scheme of finite type mapping onto $\mathbb A^{d_1|d_2}$
via  the homomorphism of rings $a_i\mapsto a_i[0]$. Next, for
an arbitrary etale $\phi: U\to\mathbb A^{d_1|d_2}$ one has,
as in \cite{KV1} (1.7.3), that 
$$\Lc^\epsilon_\phi \,\,\, = \,\,\, \Lc^\epsilon_n(\mathbb A^{d_1|d_2}) 
\times_{\mathbb A^{d_1|d_2}} U.
$$  
\end{rem}


\section{Factorization structure on $\Lc X$.}

\subsection{Reminder on $\Lc_{C^I}^0X$ and $\Lc_{C^I}X$.}\label{LCX-reminder-subsection}
Let us extend the construction of the global formal loops space from
\cite{KV1} to the case of targets belonging to the super category.
Let $C$ be a (purely even) smooth algebraic curve
 and $X$ be a smooth algebraic super-variety.
Let $\Fsetb$ be the category of nonempty finite sets
and their surjections.
Let $I$ belong to $\Fsetb$.
Let $S$ be a super-scheme and $c_I: S\to C^I$
be a morphism, so $c_I = (c_i: S\to C)_{i\in I}$. Let
$\Gamma_i\subset S\times C$ be the graph of $c_i$. Let
$\Gamma = \bigcup_{i\in I}\Gamma_i$ be the union. 
We denote by $\widehat{\Oc}_\Gamma$
the completion of $\Oc_{S\times C}$ along $\Gamma$, and by
$\Kc_\Gamma$ the localization $\widehat{\Oc}_\Gamma [r^{-1}]$, where
$r$ is a local equation of $\Gamma$ in $S\times C$. Finally,
let $\Gamma_{\red} = \Gamma\cap(S_\red \times C)$, and
$\Kc_\Gamma^\sqr \subset \Kc_\Gamma$ be the subsheaf formed
by sections whose restriction to $S_\red \times C$ lies in
$\widehat{\Oc}_{\Gamma_\red}$. 

\begin{lem} Let $\ul{\Gamma}$ be the underlying topological space of
the superscheme $\Gamma$. Then $\widehat\Oc_\Gamma$ and
$\Kc_\Gamma^\sqr$ are sheaves of local supercommutative algebras on $\ul{\Gamma}$,
so $(\ul{\Gamma}, \widehat\Oc_\Gamma)$ and $(\ul{\Gamma}, \Kc_\Gamma^\sqr)$
are super-spaces. 
\end{lem}

The proof is similar to that of Lemma \ref{nil-local-ring-lemma}. \qed

 \vskip .3cm

Consider the functor 
\be\label{lambda-functor-equation}
\lambda_{X, C^I}: S\mapsto
\biggl\{ (c_I, \phi);\,\,
c_I: S\to C^I,\,\, \phi\in \Hom_{\Sspb}\bigl((\ul{\Gamma}, \Kc_\Gamma^\sqr), \,
(\ul{X}, \Oc_X)\bigr) \biggr\}.
\ee
We define the functor $\lambda_{X, C^I}^0$ is a similar way, 
with $\widehat{\Oc}_\Gamma$ instead of $\Kc_\Gamma^\sqr$. 

\vskip .2cm

We denote by $\gen$ the Lie algebra $\Der\, \CC[[t]]$ and by
$K$ the group scheme
$$\Aut\,\CC\lsb t\rsb\,\,\, =\,\,\,
\Spec\bigl(\CC[a_1^{-1},a_1,a_2,a_3,...]\bigr). $$
So for a ring $R$ an $R$-point of $K$ is a formal change of
coordinates
$$t\mapsto a_1 t + a_2 t^2 + ..., \quad
a_1\in R^\times, \quad a_i\in R,\,\,\, i\geqslant 2.$$
The Lie algebra $\gen$ and the group scheme $K$ form a {\em Harish-Chandra pair}, see
\cite{BD1} (2.9.7).
By an action of $(\gen, K)$ on an ind-scheme $Y$ we mean an action of $K$
by automorphisms and an action of $\gen$ by derivations (infinitesimal
automorphisms) which are compatible.

Let $C$ be as before and $\widehat {C}\to C$ be the scheme whose points are
pairs $(c, t_c)$ where $c$ is a point of $C$ and $t_c$ is a formal coordinate
near $c$. The Harish-Chandra pair $(\gen, K)$ acts on $\widehat{C}$ with the
action of $K$ preserving the projection $\widehat{C}\to C$ and the action
of the element $d/dt$ of $\gen$ defining an integrable
connection on $\widehat{C}$ along $C$.

\begin{prop}\label{LCX-functor-prop}
(a) The functor $\lambda_{X, C^I}$ is represented by a super-ind-scheme 
$\Lc_{C^I}X$ over $C^I$,  
and $\lambda_{X, C^I}^0$ by a super-subscheme 
$\Lc_{C^I}^0X$ over $C^I$.  

(b)
If $I=\{1\}$, the ind-scheme $\Lc_CX$ and the scheme
$\Lc_C^0X$ are obtained by the principal bundle construction
of Gelfand-Kazhdan,  i.e.,
$$\Lc_CX=\Lc X\times_K\widehat C,
\quad
\Lc_C^0X=\Lc^0 X\times_K\widehat C.
$$
 \end{prop}
 
 \noindent {\sl Proof:} This is quite similar to \cite{KV1} (2.3-7), so we indicate
 the main steps. First, 
 we consider the case when $X=\mathbb A^1$ with coordinate $t$. As in 
 \cite{KV1}, (2.7) 
 we see that  for $S=\Spec(R)$ an affine superscheme,  a morphism $c_I: S\to C^I$
 is given by an $I$-tuple of elements $(b_i\in R)_{i\in I}$. Then
 $\Gamma_i$ is given by the equation $t=b_i$ and $\Gamma$ is given by
 $\prod_{i\in I} (t-b_i)=0$, so the completion of $\Oc_{S\times C}$ along $\Gamma$
 is described explicitly by
 $$H^0(\Gamma, \widehat{\Oc}_\Gamma)\,\,=
 \,\,\pro_{n\geqslant 0} \,\, R[t]\bigl/\prod_{i\in I} (t-b_i)^{n+1},$$
 which is then identified with 
  the set of formal series 
 \be\label{O-hat-Gamma-equation}
 \sum_{l\geqslant 0} a_l(t) \prod_{i\in I} (t-b_i)^l, \quad a_l(t)\in R[t], \,\deg(a_l)< |I|.
 \ee 
 Similarly, $H^0(\Gamma, \Kc_\Gamma)$ is identified with the set of series
 \be\label{K-Gamma-equation}
 \sum_{l\gg -\infty} a_l(t) \prod_{i\in I} (t-b_i)^l, \quad a_l(t)\in R[t], \,\deg(a_l)< |I|.
 \ee 
   The subring
 $H^0(\Gamma, \Kc_\Gamma^\sqr)$ is specified by the condition that the coefficients
 of $a_l(t)$, $l<0$, are nilpotent in $R$. 
   
  Therefore, if $X=\mathbb A^{d_1|d_2}$
   then a morphism
 $\phi$ as in \eqref{lambda-functor-equation} is just given by specifying,
 for each $l\in\ZZ$, a vector-valued polynomial $a_l^{(\phi)}(t) \in R[t]\otimes \CC^{d_1|d_2}$
 with the condition each component of each coefficient has even parity, and the 
 components of the coefficients of $a_l^{(\phi)}$ with $l<0$, are nilpotent.
 This describes $\Lc_{C^I}\mathbb A^{d_1|d_2}$ explicitly, in terms of the polynomial
 and power series rings in these  components considered as independent
 variables, as in \cite{KV1} (2.7.2). 
 Similarly for $\Lc^0_{C^I}\mathbb A^{d_1|d_2}$.
 
 Next, if $X$ is an affine super-scheme of finite type, then
 we realize $X$ as  a closed sub-superscheme of  some $\mathbb A^{d_1|d|2}$
and then realize $\Lc_{C^I}X$ inside  $\Lc_{C^I}\mathbb A^{d_1|d_2}$
by imposing the equations of $X$ identically on $d_1+d_2$-tuple
of indeterminate series \eqref{K-Gamma-equation}. Similarly for 
$\Lc^0_{C^I}\mathbb A^{d_1|d_2}$.

To treat the case of an arbitrary super-scheme of finite type, we prove
the analog of the gluing property of the functors $\lambda_{X, C^I}$ and
$\lambda^0_{X, C^I}$ as in \cite{KV1}, Proposition 2.6.1. This analog follows
directly from the definition of the functors in terms of morphisms of
superspaces as in  \eqref{lambda-functor-equation}. 

Finally, we pass from the case $C=\mathbb A^1$ to the case of an arbitrary smooth
curve by using \'etale local coordinates on $C$. This proves part (a) of the 
proposition. 

\vskip .2cm

To prove part (b), notice that for $C=\mathbb A^1$, the  choice of a
coordinate $t$ on $C$
gives a section $C\to\widehat{C}$ and thus a splitting of the Gelfand-Kazhdan 
construction, identifying,
say $\Lc X\times_K \widehat{C}$, with $\Lc X\times C$. In the presense of
such identification, the identification of $\Lc_C X$ with $\Lc X\times C$
is immediate for $X=\mathbb A^{d_1|d_2}$ and thus for $X$ closed in
$\mathbb A^{d_1|d_2}$ from the explicit construction above (the polynomials
$a_l$ will have degree 0). As the statement is local, the canonical identification
for any affine $X$ that this produces, entails an identification for any
$X$ of finite type. The case of an arbitrary $C$ can be treated
by working locally on $C$. So we can assume that $C$ has an \'etale coordinate $t$
which again splits the Gelfand-Kazhdan construction and the argument is similar. 
\qed


\subsection {Factorization structure.}

The category $\Fsetb$ has a final object $\{1\}$ (a one-point set)
and a monoidal structure $\sqcup$ (disjoint union)
but no unit object for $\sqcup$.
Let $\Fsetall$ be the category of all finite sets and all maps.
This is a monoidal category with the unit object $\varnothing$.

If $p: J\to I$ and $p': J'\to I'$
are two morphisms of $\Fsetall$, we denote
their disjoint union by
$$p\sqcup p':\,\,\, I\sqcup I' \longrightarrow  J\sqcup J'.$$
 
\vskip .2cm

Let $C$ be any super-scheme of finite type. 
For every morphism $p: J\to I$ in $\Fsetall$
we denote by  $C^p$ the open subset in $C^J$
consisting of the $J$-tuples $(c_j)$
such that $c_j\neq c_{j'}$ for $p(j)\neq p(j')$.
We will write $p_J$, or simply $J$, for the unique map $J\to\{1\}$.
Notice that $C^{p_J}=C^J$ so the two notations are compatible.
We will also write $1_J:J\to J$ for the identity.

 \vskip .2cm
 
Let $K\buildrel q\over\lra J\buildrel p\over\lra I$ be a composable pair of
morphisms of $\Fsetall$. We have the {\em diagonal map}
$$\Delta_{p,q}: \, C^p\,\, \to\,\,  C^{pq}, \quad (c_j)\,\, \mapsto\,\, (c_{q(k)}).$$
If $q$ is surjective, then $\Delta_{p,q}$ is a closed embedding. 
 We also have the {\em off-diagonal map}
$$j_{p,q}:\, C^{q}\,\,\, \to\,\,\,  C^{pq},\quad (c_k)\,\, \mapsto \,\,(c_k),$$
which is always an open embedding.
For each $p,p'$ we have also the map
$$i_{p,p'}:\,\,C^{p\sqcup p'}\to C^{p}\times C^{p'},\quad (c_k)\mapsto(c_k),$$
which is also an open embedding.
The above maps fit into the following commutative diagrams,
existing for any composable triple
$L\buildrel r\over\lra K\buildrel q\over \lra J\buildrel p\over
\lra I$ of morphisms of $\Fsetall$:

\be\label{diagonal-maps-equation}
 \begin{aligned}
\begin{matrix}  C^{p}&\buildrel\Delta_{p,q}\over\lra& C^{pq}&\cr
\Delta_{p, qr}&\searrow&\big\downarrow&\Delta_{pq,r}\cr
&&C^{pqr}&\end{matrix}, 
\\
\begin{matrix} C^{r}& \buildrel j_{q,r}\over\lra&C^{qr}&\cr
j_{pq,r}&\searrow&\big\downarrow&j_{p, qr}\cr
&&C^{pqr}&\end{matrix},
 \\
\begin{matrix} & C^{q}&\buildrel j_{p,q}\over\lra& C^{pq}&\cr
\Delta_{q,r}&\big\downarrow&&\big\downarrow&\Delta_{pq,r}\cr
&C^{qr}&\buildrel j_{p, qr}\over\lra&C^{pqr}&\end{matrix}.
\end{aligned}
\ee

\begin{Defi}\label{factorization-semigroup-definition}
Let $Y_C\to C$ be a  super-ind-scheme formally smooth over $C$, 
equipped with an integrable connection along $C$.
A {\em factorization semigroup} on $Y_C$ is a system consisting of :

(a)  For any morphism $p$ of $\Fsetb$, 
a super-ind-scheme $\rho_p:Y_p\to C^{p}$ formally smooth
over $C^{(p)}$,
equipped with integrable connections along $C^{p}$,
so that $Y_{\{1\}}=Y_C$,

(b) for any composable pair $p,q$
in $\Fsetb$, morphisms of
relative super-ind-schemes with connections
$$\varkappa_{p,q}:\,\, \Delta_{p,q}^* (Y_{pq})\,\,\, \to\,\,\, Y_p,
\quad \kappa_{p,q}: \,\,j_{p,q}^* (Y_{pq})\,\,\,\to\,\,\, Y_q$$
which are isomorphisms and
satisfy the compatibility conditions lifting 
\eqref{diagonal-maps-equation}:
$$\varkappa_{p, qr}\,\,\, = \,\,\,\varkappa_{p,q}\circ\Delta_{p,q}^*(\varkappa_{pq,r}): \,\,
\Delta_{p, qr}^* (Y_{pqr})\,\,\,\to\,\,\, Y_p,
$$
$$\kappa_{pq,r}\,\,\,=\,\,\,
\kappa_{q,r}\circ j_{q,r}^*(\kappa_{p, qr}): \,\,j_{pq,r}^*(Y_{pqr})\,\,\,\to\,\,\, Y_r,
$$
$$
\kappa_{p,q}\circ j_{p,q}^*(\varkappa_{pq,r})\,\,\,=\,\,\, \varkappa_{q,r}\circ\Delta_{q,r}^*
(\kappa_{p, qr}):\,\, j_{p,q}^*\Delta_{pq,r}^*(Y_{pqr})\,\,\,=\,\,\,
\Delta_{q,r}^* j_{p, qr}^*(Y_{pqr})\,\,\, \to\,\,\, Y_q,
$$
(c)
for any pair $p,p'$ in $\Fsetb$, isomorphisms
$$\sigma_{p,p'}:i_{p,p'}^*(Y_p\times Y_{p'})\to Y_{p\sqcup p'}.$$
\end{Defi}

\begin{Defi}\label{cocommutative-semigroup-definition}
A factorization semigroup $(\rho_p:Y_p\to C^p)$ is said to be {\em cocommutative}
if, for any $J, J'$  the maps $\varkappa$, $\kappa$ factor through a morphism
of $C^{J\sqcup J'}$-schemes
$Y_{J\sqcup J'}\to Y_{J}\times Y_{J'}$. Here $Y_J = Y_{p_J}$.
\end{Defi}

\begin{exa}
The collection $(C^p)$ forms a cocommutative factorization semigroup which
we call the {\em unit semigroup}.
\end{exa}

In the remainder of this subsection we will assume that  $C$ is a purely
even smooth algebraic curve.

\begin{rems}\label{semigroups-vs-monoids-remarks} {\bf Semigroups versus monoids.}  
(a) 
The definition \eqref{factorization-semigroup-definition} is equivalent to
 \cite{KV1} (2.2.1).
Indeed, given a system $(Y_p)$ as before,
we define $Y_I = Y_{p_I}$.
Then the $Y_I$ satisfy the conditions of {\it loc. cit}.
Conversely, given $(Y_I)$ as in {\it loc. cit.} and $p: J\to I$ a surjection,
we define $Y_p= j_{p_I,p}^*(Y_J)$.
Then the $Y_p$ satisfy the conditions of \eqref{factorization-semigroup-definition}.
The reason for the definition chosen here is that it allows one to easily treat
higher compatibility conditions, which become necessary when
dealing with factorizing line bundles, factorizing
gerbes etc. This will be important in the subsequent paper. 

\vskip .2cm

(b) In this paper we changed the terminology of \cite{KV1}
by calling factorization semigroups what was there called
{\em factorization monoids}. Indeed, it is more natural,
following \cite{BD1} (3.10.16), to
reserve the term
``factorization monoid" to mean
a similar structure, but with $Y_p$ defined for any
morphism $p$ in $\Fsetall$, the morphisms $\kappa_{p,q}$ and
$\sigma_{p,p'}$ being always isomorphisms, and $\varkappa_{p,q}$
being an isomorphism for surjective $q$. A factorization
monoid $(Y_p)$ possesses a {\em unit section} which is a
collection of sections $(e_p: C^p\to Y_p)$, $p: J\to I$, defined as follows.
Take $q: \varnothing\to J$, then $C^{pq}=\{\bullet\}$, and the
analog of the axiom (c) implies that $Y_{pq}=\{\bullet\}$ as well. Thus
$\Delta_{p,q}^*(Y_{pq})=C^p$, and  $\varkappa_{p,q}$
is a morphism from $C^p$ to $Y_p$. We define $e_p$ to be this morphism. 
It then follows, in particular, that $(e_p: C^p\to Y_p)$
is a morphism of factorization semigroups. It also follows that
 for any local section $s$ of $Y_C\to C$, the product
$y_{\{1\}}\times s$ extends to a section of $Y_{\{1,2\}}$
(via $\kappa$, $\sigma$)
whose restriction to the diagonal is identified with $s$
(via $\varkappa$).




\vskip .2cm

(c) One can compare our concept of a factorization monoid/semigroup
with that of a chiral monoid/semigroup as introduced in \cite{BD1} (3.10.16).
The latter objects live on symmetric powers of $C$, not Cartesian powers.
In addition, the authors impose a condition which (translated
into the Cartesian power language) means that
the closure in
$Y_I$ of the complement to the preimage of the discriminant divisor in
$C^I$ equals $Y_I$.

\vskip .2cm

 (d) The map in Definition \ref{cocommutative-semigroup-definition}
 goes in the direction opposite to the map
in \cite{BD1} (3.10.16)
in the axioms of commutative chiral monoids.

\vskip .2cm

(e) The integrable connection of a factorization monoid can be recovered from
the other axioms as follows.
Assume that $p=p_{\{1\}}$, so $C^p=C$ and $Y_p=Y_C$. Let us show how to recover
the connection on $Y_p$ in this case. The general case is similar. We use the second description
of integrable connections in Proposition \ref{int-connections-conditions-prop}.
Set $J=\{1,2\}$, so $C^{[J]}$ is the formal neighborhood of the diagonal
in $C^2$, and let
$q_1, q_2:C^{[J]}\to C$ be the coordinate projections.
 We will construct an isomorphism of super-ind-$C^{[J]}$-schemes
$q_1^*(Y_C)\to q_2^*(Y_C)$ which restricts to the identity
of $Y_C$ over the diagonal $C\subset C^J$.
By definition of the unit, the maps
$\Id\times y_{\{2\}}$, $y_{\{1\}}\times\Id$
yield isomorphisms
$q_1^*(Y_{C})\to (Y_J)|_{C^{[J]}}$, $q_2^*(Y_{C})\to (Y_J)|_{C^{[J]}}$
which restrict to the identity over the diagonal.
This gives a connection.
Further, taking $I = \{1,2,3\}$
and using the unit property gives at once the integrability, see
Section 3.4.7 of \cite{BD1}.
\end{rems}

Now, let $X$ be a smooth algebraic super-variety.
Recall that  $C$ is a purely
even smooth algebraic curve.
Given a morphism $p:J\to I$ in $\Fsetb$, 
we denote by $\Lc_pX$, $\Lc_p^0X$ the restrictions of
$\Lc_{C^J}X$, $\Lc_{C^J}^0X$ to the subscheme $C^{p}$ of $C^J$.
We have the morphisms
\be 
X\buildrel \pi_p\over\lla \Lc^0_pX\buildrel i_p\over\longrightarrow\Lc_pX
\buildrel\rho_p\over\longrightarrow C^{p}.
\ee

\begin{prop}\label{LX-factorization-prop}
The systems $(\Lc^0_pX)$, $(\Lc_pX)$ 
are structures of factorization semigroups
  on $\Lc^0_CX$, $\Lc_CX$. Further, the factorization semigroup
$(\Lc^0_C X)$ is cocommutative. 
\end{prop}

\noindent {\sl Proof:} In the case of purely even $X$,
the factorization structure was given
in \cite{KV1} (2.3.3) and established at the level of the functors
$\lambda_{X, C^I}$ and $\lambda^0_{X, C^I}$ represented by
$\Lc_{C^I}X$ and $\Lc^0_{C^I}X$. This argument extends verbatim to
the case when $X$ is a smooth algebraic super-variety. 

\vskip .2cm

However, the integrable connections along $C^p$
were not given in \cite{KV1}. Since 
the factorization semigroups $(\Lc^0_pX)$, $(\Lc_pX)$ have no units,
in the sense of Remarks \ref{semigroups-vs-monoids-remarks}, these  connections 
have to be defined separately. Here we supply the definition.
We use the  formulation of  an integrable connection as a datum of type (1) in
Proposition \ref{int-connections-conditions-prop}. We will construct the
connection on $Y_I = Y_{p_I}\to C^I$, and the case of arbitrary $Y_p$, $p: J\to I$ will
follow by restriction to an open subset $C^p\subset C^J$. 

So let $c_I, c'_I: S\to C^I$ be two infinitely near maps, with components
$c_i, c'_i: S\to C$, $i\in I$. Constructing the data in
Proposition \ref{int-connections-conditions-prop}(1),
we will explain how to canonically identify the
pullback $c_I^* \Lc_{C^I}X$ with $c'_I{}^*\Lc_{C^I}X$, and similarly for the pullbacks
$\Lc^0_{C^I}X$. Indeed, 
for each $i$ we have that $c_i$ and $c'_i$
are infinitely near. Let $\Gamma, \Gamma'\subset S\times C$ be the graph unions for
$c_I$ and $c'_I$. Note that the underlying topological spaces of $\Gamma$ and $\Gamma'$
are the same. By definition
\eqref{lambda-functor-equation},
 a morphism from  $S$ to $\Lc_{C^I}X$  covering $c_I$, is the same as  a morphism 
 of superspaces $(\ul{\Gamma}, \Kc_\Gamma^\sqr)\to (\ul{X}, \Oc_X)$. Similarly,  a morphism
 $S\to \Lc^0_{C^I}X$ covering $c_I$, is the same a morphism
 $(\ul{\Gamma}, \widehat\Oc_\Gamma)\to (\ul{X},\Oc_X)$. Therefore, in order to
 identify the pullbacks, it is enough to prove the following:
 
  \begin{lem}\label{identification-of-K-lemma}
 In the situation described, we have a canonical identification of sheaves 
 on
 $\ul{\Gamma}=\ul{\Gamma'}$:
$$\widehat\Oc_\Gamma \simeq \widehat\Oc_{\Gamma'},\quad 
 \Kc_{\Gamma }^\sqr\simeq\Kc_{\Gamma'}^\sqr.$$
 
\end{lem}

\noindent{\sl Proof:} This statement is local on $C$. 
Choosing an \'etale coordinate on $C$ we reduce to the case $C=\mathbb A^1$,
sofor each $i$ we can see $c_i$, $c'_i$ as elements of the coordinate ring 
$B:=\CC[S]$ such that
$s_i=c_1-c_2$ is a nilpotent element of $B$. 
Let $n_0$ be such that $s_i^{n_0}=0$.
Put $R=B[t]=\CC[S\times\mathbb A^1]$ and let $r_i=t-c_i$, $r'_i=t-c'_i$, so
$r_i-r'_i= c'i-c_i =-s_i$. 
The equation of $\Gamma$
in $S\times C$ is then $r=\prod_i r_i$, while the equation of $\Gamma'$
is $r'=\prod_i r'_i$. 
Then
$$\widehat \Oc_\Gamma \,\,=\,\,\pro_n \Oc_{S\times C}/(r^n), \quad 
\widehat \Oc_{\Gamma'} \,\,=\,\,\pro_n \Oc_{S\times C}/(r'{}^n).$$
 On the other hand, since $r_i-r'_i$ is nilpotent for each $i$, so is $r-r'$.
This implies that 
 the $r_1$-adic and the $r_2$-adic topologies on $\Oc_{S\times C}$ are equivalent to
each other.
This implies the first identification of the lemma.

\vskip .2cm

To prove the second identification, we recall that
$\Kc_\Gamma$ is obtained from $\widehat\Oc_\Gamma$
by inverting a local equation of $\Gamma$ which we can take to be the element $r$ above.
Similarly for $\Kc_{\Gamma'}$ and $r'$. 
Let
$$\widehat R\,\, =\,\, \pro \,R/(r^n)\,\,=\,\,\pro \, R/(r'{}^n).$$
Then $r,r'\in\widehat R$ with $s=r-r'$ nilpotent, 
and it is enough to prove  that
$$\widehat R[r^{-1}]=
\widehat R[r'{}^{-1}] .$$
To see this, let us write
$${1\over r}={1\over r'}\Bigl(1-{s\over r'}+{s^2\over r'{}^2}-\cdots\Bigr)$$
(a terminating geometric series). So $r$ is invertible in
$\widehat R[r'{}^{-1}]$. Changing the order, we see that $r'$ is invertible in
$\widehat R[r{}^{-1}]$. This implies that $\Kc_\Gamma$ is identified with
$\Kc_{\Gamma'}$. Further, the subsheaves $\Kc_\Gamma^\sqr$ and $\Kc_{\Gamma'}^\sqr$
are defined by the condition involving restriction to $S_\red\times C$, and
are therefore identified as well. 
\qed

\vskip .2cm

We finally explain why $(\Lc^0_pX)$ gives a cocommutative factorization semigroup
structure. In other words, for each nonempty finite sets $I, I'$ we construct
morphisms of $C^{I\sqcup I'}$-ind-schemes
\be 
\Lc^0_{C^{I\sqcup I'}}X \,\,\,\to\,\,\, 
\Lc^0_{C^I}X \,\times \Lc^0_{C^{I'}}X.
\ee
Indeed, let $S$ be a super-scheme and $c_{I\sqcup I'} = (c_I, c_{I'})$
be a morphism from $S$ to $C^{I\sqcup I'}$. Let $\Gamma_I$
be the union of the graph of the components $c_i: S\to C$ of $c_I$. Similarly 
$\Gamma_{I'}$, $\Gamma_{I\sqcup I'}$. Now, the second datum of a morphism
from $S$ to $\Lc_{C^{I\sqcup I'}}X$ is a morphism $\phi$ from
the formal neighborhood of $\Gamma_{I\sqcup I'}$ in $S\times C$,  to $X$.
Now, the formal neighborhoods of $\Gamma_I$ and $\Gamma_{I'}$  are each contained in that
of $\Gamma_{I\sqcup I'}$, so by restricting $\phi$ we get morphisms of
these formal neighborhoods into $X$ which, together with the $c_I$, $c_{I'}$
give morphisms $S\to \Lc_{C^I}^0X$ and $S\to \Lc^0_{C^{I'}}X$.
This finishes the proof of Proposition \ref{LX-factorization-prop}.

\begin{rem}\label{LSX=SLX-remark} 
The same proof as in Proposition \ref{LSX=SLX-proposition} implies that,
for any super-scheme $X$ of finite type, we have an isomorphism of super-ind-schemes
$$\Lc_{C^p}\Sc X=\Sc \Lc_{C^p}X.$$ Indeed, both represent the functor
\begin{equation*}
S\mapsto
\biggl\{ (c_p, \phi);\,\,
c_p: S\to C^p,\,\, \phi\in \Hom_{\Sspb}\bigl((\ul{\Gamma}, \Kc_\Gamma^\sqr[\eta]), \,
(\ul{X}, \Oc_X)\bigr) \biggr\}
\end{equation*}
where $\eta$ is an  odd generator.
\end{rem}


\subsection{Factorization of $\Lc \Sc^N X$ on super-curves.}\label{super-curves-factorization-subsection}

Fix an integer $N\geqslant 0$. 
Let now $C$ be a smooth super-curve of pure dimension $(1|N)$. For every $\CC$-point
$c\in C$ the completed local ring $\widehat\Oc_{C, c}$ is isomorphic to
$\CC\lsb t\rsb [\eta_1, ..., \eta_N]$. More generally, let $c: S\to C$ be a
point of $C$ with values in a super-scheme $S$. Denoting $\Gamma_c\subset S\times C$
the graph of $c$, we have the completion $\widehat\Oc_c$ of $\Oc_{S\times C}$ along
$\Gamma$, and we call a {\em formal coordinate system} at $c$ an isomorphism
of sheaves of topological local rings
$$\Oc_\Gamma\lsb t\rsb [\eta_1, ..., \eta_N]\lra \widehat\Oc_c.$$
As in Subsection \ref{LCX-reminder-subsection}, we have a super-scheme
$\widehat C\to C$ whose $S$-points
are  data $(c, t, \eta_1, ..., \eta_N)$ consisting of an $S$-point $c: S\to C$
and a formal coordinate system $(t,\eta_1, ..., \eta_N)$ at $c$.  We also have
the Harish-Chandra pair $(\gen_{1|N}, K_{1|N})$.  Here $\gen_{1|N}$ is the
Lie super-algebra $\Der \CC\lsb t\rsb [\eta_1, ..., \eta_N]$ while
$K_{1|N} = \Aut \,\CC \lsb t\rsb [\eta_1, ..., \eta_N]$ is the group super-scheme
whose points in a super-commutative algebra $R$ are invertible formal changes of coordinates
$$t\mapsto \sum_{i\geqslant 0; \, J} a_{i,J} t^i \eta^J, \quad
\eta_\nu\mapsto \sum_{i\geqslant 0; J} b^\nu_{i, J} t^i\eta^J.$$
Here $J$ runs over subsets  $J=\{ 1\leqslant j_1< ... < j_p\leqslant N\}$, 
the element
$a_{i,J}\in R$ is of parity $|J|$, 
the element $b^\nu_{i,J}$ is of parity $|I|+1$, and we have
$$\eta^J=\eta_{j_1}\eta_{j_2}\dots\eta_{j_p}.$$ 
It is further required that $a_{0,\varnothing}=0$, $a_{1,\varnothing}\in R_\even^*$, and
the matrix $\| b^\nu_{0, \{\mu\}}\|_{\mu,\nu=1}^N$ is invertible. This Harish-Chandra pair
acts on $\widehat C\to C$ as in the even case. As in Proposition \ref{LSX=SLX-proposition} the ind-scheme
$\Lc\Sc^N X= \Sc^N\Lc X$ represents the functor
$$S\,\,\mapsto\,\,\Hom_\Sspb\bigl((\ul{S},\Oc_S\lb t\rb^\sqr[\eta_1, ..., \eta_N]),(\ul{X},\Oc_X)\bigr),$$
where $\eta_1, ..., \eta_N$ are  odd generators. Similarly for $\Lc^0\Sc^N X$ and
$\Oc_S\lsb t\rsb [\eta_1, ..., \eta_N]$. 
The Harish-Chandra pair $(\gen_{1|N}, K_{1|N})$ also acts on the super-scheme
$\Lc^0\Sc^N X$ and the super-ind-scheme $\Lc\Sc^N X$, thus giving
a super-scheme  and a super-ind-scheme
\be\label{super-LCX}
\Lc^0_CX = \Lc^0\Sc^N X\times_{K_{1|N}} \widehat C\lra C,
\quad \Lc_C X= \Lc\Sc^N X\times_{K_{1|N}} \widehat C
\lra C
\ee
with
integrable connections along $C$. These integrable connections are given
by the action of $\partial/\partial t, \partial/\partial\eta_\nu\in\gen_{1|N}$. 

\begin{prop}\label{factorization-supercurves-prop} 
For any $N\geqslant 0$ and any smooth super-curve $C$ of
dimension $(1|N)$ there  exist factorization semigroups $(\Lc^0_p X)$, resp. 
$(\Lc_p X)$ on  the super-ind-schemes $\Lc^0_C X$, resp. $\Lc_C X$
given by \eqref{super-LCX}. 
\end{prop}

\noindent {\sl Proof:} The construction is similar to that of Subsection 
\ref{LCX-reminder-subsection}. That is, for any $I\in\Fsetb$ and any
$c_I: S\to C^I$  with components $c_i: S\to C$, we denote by $\Gamma\subset S\times C$
the union of the graphs of the $c_i$ and construct three sheaves
$\widehat\Oc_\Gamma, \Kc_\Gamma$, and $\Kc_\Gamma^\sqr$ on the underlying
topological space $\ul{\Gamma}$. Of these, $\widehat\Oc_\Gamma$ is just the
completion of $\Oc_{S\times C}$ along $\Gamma$ (so its construction
does not use the specifics of $C$ being a super-surve). Next,
the definition of $\Kc_\Gamma$ is based on the following lemma.

\begin{lem} Let $(r,\xi_1, ..., \xi_N)$ and $(r',\xi'_1, ..., \xi'_N)$
be two systems of local equations for $\Gamma$ in $S\times C$,
with $r,r'$ being even and $\xi_\nu, \xi'_\nu$ being odd. Then
$r'$ is invertible in
$\widehat\Oc_\Gamma[r^{-1}]$, and $r$ is invertible in $\widehat\Oc_\Gamma[r'{}^{-1}]$.
\end{lem} 

\noindent {\sl Proof:} follows from the nilpotency of $r-r'$, as in 
Lemma \ref{identification-of-K-lemma}.\qed

\vskip .2cm

The lemma implies that we have a well-defined sheaf $\Kc_\Gamma=\widehat\Oc_\Gamma[r^{-1}]$,
and we define $\Kc_\Gamma^\sqr\subset\Kc_\Gamma$ as in Subsection 
\ref{LCX-reminder-subsection}. After this we define the functor $\lambda_{X, C^I}$
as in \eqref{lambda-functor-equation}, using $\Kc_\Gamma^\sqr$ and similarly
for $\lambda^0_{X, C^I}$ using $\widehat\Oc_\Gamma$. The proof of representability of
these functors is completely analogous to the proof of Proposition 
\ref{LCX-functor-prop}(a). Finally the proof that these functors yield factorization semigroups is
completely analogous to the proof of Proposition \ref{LX-factorization-prop}.
We leave the remaining details to the reader. \qed
 
\begin{rem} Factorization semigroups on $(1|N)$-dimensional
super-curves are non-linear analogs of  $N_W=N$ SUSY vertex algebras
as defined in \cite{Heluani}. More precisely, recall that the categories of factorization algebras
and chiral algebras on a curve are equivalent \cite{BD1}. 
One can define a category of factorization algebras on a given $(1|N)$-dimensional
super-curve $C$ which is equivalent to the category of chiral algebras on $C$ considered in
\cite{Heluani}.
Further, factorization semigroups on $C$ yield natural examples of factorization algebras on $C$,
and $N_W=N$ SUSY vertex algebras yield chiral algebras on $C$ according to \cite{Heluani}.
In particular, Proposition \ref{factorization-supercurves-prop}
provides a geometric reason for the observation of \cite{Ben-Zvi-Heluani-Szczesny}
that  $\Omega^{\on{ch}}_X$, the chiral de Rham complex of any manifold, is a sheaf of
$N_K=1$ SUSY vertex algebras. Indeed, $\Omega^{\on{ch}}_X$ can be seen
as a sheaf of chiral differential operators on $\Sc X$ and can be
recovered from 
$\Lc\Sc X$ and its $(1|1)$-dimensional factorization structure.

\end{rem}


\section{The transgression.}

\subsection{Definition of the transgression.}

Recall from Subsection \ref{functions-forms-superschemes-subsection}
 that for every super-ind-scheme $Y$ we have
a sheaf $\Omega^m_Y$ on the topological space $\ul{Y}$.
In particular, if $Y=\Lc X$, then $\ul{Y}=\ul{\Lc^0X}$.
We define
$$\Omega^m_{\Lc X|\Lc^0X}=\Ker\{\Omega^m_{\Lc X}\to\Omega^m_{\Lc^0X}\}.$$
In particular, for $m=0$ we write
$$\Oc_{\Lc X|\Lc^0X}=\Omega^0_{\Lc X|\Lc^0X}.$$

Let $R = \pro{}_{\a\in A} \, R_\a$ be a super-commutative
pro-algebra, or, what is the same,
a topological super-algebra represented as a filtering projective limit of
discrete super-commutative
algebras $R_\a$. The ring of Laurent series with coefficients
in $R$ is defined by
\be 
R\lb t\rb \,\,\, = \,\,\, \pro_{\a\in A} \,\, R_\a \lb t\rb \,\,\, = \,\,\,
\biggl\{ \sum_{n=-\infty}^\infty a_n t^n \biggl| \, a_n\in R, \,\, \lim_{n\to-\infty} a_n
=0\biggr\}. 
\ee

As in \cite{KV2} (6.2), we have the evaluation map
which is a morphism of ringed spaces
\be 
\ev: \,\,\bigl( \ul{\Lc^0X}, \Oc_{\Lc X}\lb t\rb\bigr) \,\,\,\lra\,\,\,
(\ul{X}, \Oc_X).
\ee
Its underlying morphism of topological spaces is
$$\ev_\flat = \pi: \,\, \ul{\Lc^0X}\lra \ul{X}.$$
In terms of the identification of $\pi_*\Oc_{\Lc X}$ given in
Proposition \ref{borisov-description-LX-prop}(b), the morphism of sheaves of rings
corresponding to $\ev$ is
\be 
\ev^\sharp: \,\, \pi^{-1}\Oc_X \,\,\lra\,\, \Oc_{\Lc X}\lb t\rb,\quad 
\ev^\sharp (a) \,\,\, = \,\,\, \sum_{n=-\infty}^\infty a[n] t^n.
\ee

\begin{rem}
Let $S^1$ be the unit circle $|t|=1$ 
in the complex plane, $M$ be a complex analytic manifold, 
and $LM = C^\infty(S^1, M)$ be the space
of $C^\infty$-maps from $S^1$ to $M$. The map $\ev$ is the algebraic
analog of the canonical map
$S^1\times LM \,\,\to\,\, M.$
\end{rem}

Let us now construct, similarly to \cite{KV3} (1.3), the {\em transgression map}
\be 
\tau: \,\, \Omega^m_X \to \pi_* \Omega^{m-1}_{\Lc X|\Lc^0X},
 \ee
compatible with the differential.  
For a topological super-commutative algebra $R$ as above we have the {\em residue homomorphism}
$$\Res:\Omega^m(R\lb t\rb )\to\Omega^{m-1}(R),$$ 
see e.g., \cite{KV3}, (1.3.4) for the commutative case, the super-commutative case is given by the same
formulas.
Now, the map $\tau$ is the composition of
$$\Res:\pi_*\Omega^m(\Oc_{\Lc X}\lb t\rb )\to\pi_*\Omega^{m-1}(\Oc_{\Lc X})=
\pi_*\Omega^{m-1}_{\Lc X},$$
and the pull-back with respect to the evaluation map
$$\ev^*:\Omega^m_X=\Omega^m(\Oc_X)\to\pi_*\Omega^m(\Oc_{\Lc X}\lb t\rb ).$$

We now assume that $C$ is a purely even smooth algebraic curve and use the
factorization semigroups $(\Lc_pX)$ and $(\Lc^0_pX)$ from Subsection
\ref{LCX-reminder-subsection}.

\begin{Defi}\label{additive-forms-definition}
Let $\xi\in\Omega^m_{\Lc X|\Lc^0X}$ be a globally defined
$m$-form vanishing
on $\Lc^0X$. We say that $\xi$ is {\em additive}, if, first of all, it is
$(\gen, K)$-invariant 
and  so gives rise
to a relative $m$-form $\xi_C \in\Omega^m_{(\Lc_C X|\Lc^0_CX)/C}$.
Second, we require that 
there exists
a family $\xi_p$ of relative forms on $\Lc_pX$ over 
$C^p$, vanishing along $\Lc^0_pX$ and satisfying the conditions:

(a) For $p=\{1\}$ (the identity map of a 1-element set), we have $\xi_p = \xi_C$. 

(b) For any two composable morphisms $p,q$ of $\Fsetb$ we have
$$\Delta_{p,q}^*(\xi_{pq})\,\,\, = \,\,\,
\varkappa_{p,q}^*(\xi_p), \quad
j_{p,q}^*(\xi_{pq}) \,\,\, = \,\,\,
\kappa_{p,q}^*(\xi_q).$$

(c) For any two morphisms $p, p'$ of $\Fsetb$ we have
$$i_{p,p'}^*(\xi_p \boxplus \xi_{p'}) \,\,\, = \,\,\,
\sigma_{p, p'}^*(\xi_{p\sqcup p'}),$$
where $\boxplus$ means the differential form on the Cartesian product 
obtained by adding the pullbacks of two forms from the factors.
\end{Defi}

Note that the forms $\xi_p$, if they exist, are uniquely defined by 
the conditions above. 
We denote by $\Ac dd^m(X)$ the space of additive $m$-forms
on $\Lc X$, and by $\Ac dd^m _X$ the sheaf
\be 
U\mapsto  \Ac dd ^m(U)
\ee
on  the Zariski topology of $X$. Note that the differential of an additive form
is again additive, so we have the de Rham complex $\Ac dd^\bullet_X$
of additive forms. 
For $m=0$ we will speak of {\em additive functions}
and denote by $\Ac dd_X = \Ac dd_X^0$ the sheaf of such functions.

\begin{prop}\label{transgression-additive-proposition}
For any $m$-form $\eta\in\Omega^m(U)$, the
$(m-1)$-form $\tau(\eta)$ is additive. We have therefore a morphism of
complexes of sheaves
$$\tau: \,\, \Omega^\bullet_X \,\,\lra \,\, \Ac dd^{\bullet-1}_X.$$
\end{prop}

\noindent{\sl Proof :}
First of all, the fact that $\tau(\eta)$ is $(\gen, K)$-invariant, is clear,
as integration of differential forms is an invariant procedure. Further, 
generalizing \cite{KV3} (1.6), 
we construct the ``global" version of the transgression
map
\be\label{tau-I-equation}
\tau_I: \,\, \Omega^m_X \,\,\,\lra \,\,\,
\pi_{I*}\bigl(\Omega^{m-1}_{(\Lc_{C^I}X/C^I)|\Lc^0_{C^I}X} \bigr). 
\ee
Here the sheaf in the right hand side consists of relative $(m-1)$ forms on $\Lc_{C^I}X$
over $C^I$, vanishing along $\Lc^0_{C^I}X$, and
$\pi_I$ is the canonical projection 
$$\pi_I: \,\, \Lc^0_{C^I}X \,\,\, \lra\,\,\, X.$$

Let $\eta$ be a  local section of $\Omega^m_X$. Restricting $X$ if
necessary, we can assume that $\eta$ is a global section. To define $\tau_I(\eta)$,
we need to define, for each super-scheme $S$ and each morphism $h: S\to \Lc_{C^I}X$,
an $(m-1)$-form $h^*\tau_I(\eta)$ on $S$ in a compatible way. 
Let $h$ correspond to a datum $(c_I, \phi)$ with respect to the
graph subscheme $\Gamma$, as in
\eqref{lambda-functor-equation}. 
We then get a section
$$\phi^*\eta \,\,\,\in\,\,\,
H^0(\Gamma, \Kc_\Gamma^\sqr \otimes\Omega^m_{S\times C}).$$
Let $q: S\times C\to S$ be the projection. Then
$$\Omega^1_{S\times C} \,\,\, = \,\,\, \Omega^1_{S\times C/S} \oplus
q^*\Omega^1_S,$$
which implies that
$$\Omega^m_{S\times C} \,\,\, = \,\,\, \bigoplus_{i+j=m} \Omega^i_{S\times C/S}
\otimes q^*\Omega^j_S.$$
Let
\be\label{projection-v-equation}
v:\,\, \Omega^m_{S\times C}\,\,\, \lra \,\,\,
\Omega^1_{S\times C/S}\otimes q^*\Omega^{m-1}_S
\ee
be the projection to the summand with $i=1, j=m-1$. 

Let us denote the projection $\Gamma\to S$ by the same letter
$q$. Our statement now follows from the next lemma.

\begin{lem}\label{Res-Gamma/S-lemma}
 For each super-scheme $S$ and each morphism of super-schemes
$c_I: S\to C^I$, $c_I=(c_i: S\to C)_{i\in I}$,
there is a morphism
$$\on{Res}_{\Gamma/S}: q_*\bigl(\Kc_\Gamma\otimes\Omega^1_{S\times C/S}\bigr)
\lra\Oc_S$$
of sheaves on $S$, and these morphisms satisfy the following properties:

(a) Compatibility with base change for any morphism of super-schemes $S'\to S$.

(b) Additivity: Let $I',I''$ be two non-empty finite sets, and $c_{I'}: S\to C^{I'}$
and $c_{I''}: S\to C^{I''}$ be two morphisms whose graph unions $\Gamma', \Gamma''$
are disjoint. Denote $I=I'\sqcup I''$ and let $c_I=(c_{I'}, c_{I''}): S\to C^I$
be the combined morphism whose graph union is $\Gamma=\Gamma'\sqcup \Gamma''$.
Then, with respect to the identification
$$q_*\bigl(\Kc_\Gamma\otimes\Omega^1_{S\times C/S}\bigr)\,\,=\,\,
q'_*\bigl(\Kc_{\Gamma'}\otimes\Omega^1_{S\times C/S}\bigr)\,\oplus\,
q''_*\bigl(\Kc_{\Gamma''}\otimes\Omega^1_{S\times C/S}\bigr),$$
we have
$$\on{Res}_{\Gamma/S}(\omega'\oplus\omega'')\,\,=\,\,\on{Res}_{\Gamma'/S}(\omega')\,+\,
\on{Res}_{\Gamma''/S}(\omega'').$$

(c) Normalization: If $|I|=1$, so that $q: \Gamma\to S$ is an isomorphism, and $t$ is a local
equation of $\Gamma$ in $S\times C$, then 
$$\on{Res}_{\Gamma/S}\,\, \biggl(\sum_{n\gg-\infty}^\infty u_n t^n dt\biggr)\,\,=\,\, u_{-1}.$$
\end{lem}

Indeed, suppose we know the lemma. We then define the form 
$$h^*\tau_I(\eta)\,\,=\,\,(\on{Res}_{\Gamma/S}\otimes \Sigma)(v(\phi^*\eta))$$
on $S$ for each $S$ and each $h: S\to \Lc_{C^I}X$. Here
$\Sigma: q_*q^*\Omega^{m-1}_S\to\Omega^{m-1}_S$ is the ``trace" morphism
(summation over the fibers). 
By part (a) of the lemma, this means that we have the
form $\tau_I(\eta)$ on $\Lc_{C^I}X$, as in \eqref{tau-I-equation}. 
Let $p: J\to I$ be a morphism of $\Fsetb$. We define the $(m-1)$-form
$\tau_p(\eta)$ on $\Lc_p X$ to be the restriction of $\tau_J(\eta)$
to the open part $\Lc_p X \subset \Lc_{C^J}X$. After that, the condition (a)
of Definition \ref{additive-forms-definition} follows from part (c) of the lemma,
condition (b) follows from the definition of $\tau_p\eta$ as the restriction,
while condition (c) follows from part (b) of the lemma. 
So $\tau(\eta)$ is indeed
an additive form. The fact that $\tau$ is a morphism of complexes, i.e., 
$$\tau_p(\eta+\eta') \,\, = \,\,\tau_p(\eta) + \tau_p(\eta'), 
\quad \tau_p(d\eta) = d\tau_p(\eta),$$
follows from the 
corresponding properties of $\tau$ and from the fact that the form 
$\tau_p(\eta)$ is defined by $\tau(\eta)$ and by the conditions (a)-(c)
of Definition  \ref{additive-forms-definition} uniquely (the question being
only its existence).

\vskip .3cm

\noindent {\sl Proof of Lemma  \ref{Res-Gamma/S-lemma}, } \underbar{Step 1: $S$ is a scheme:}
In this case the construction of $\on{Res}_{\Gamma/S}$ is deduced from the
Grothendieck duality theory \cite{Conrad}, as described in \cite{KV3} (1.6).
That is, we have the principal part morphism 
$$P: \,\, \Kc_\Gamma \otimes\Omega^{1}_{S\times C/S} \,\,\,\lra\,\,\,
\Kc_\Gamma \otimes\Omega^{1}_{S\times C/S}\bigl/ \widehat\Oc_\Gamma
\otimes\Omega^{1}_{S\times C/S}\,\,=\,\,
\underline{H}^1_\Gamma(\Omega^1_{S\times C/S}),
$$
which we compose with the trace map of the Grothendieck duality
$$\on{tr}_{\Gamma/S}: q_* \underline{H}^1_\Gamma(\Omega^1_{S\times C/S})
\lra R^1q_*(\Omega^1_{S\times C/S})
\lra \Oc_S.$$
Now, compatibility of the trace map with arbitrary base change for schemes
was established in \cite{Conrad} (1.1.3), by reduction to the case of Noetherian
base ($S$ in our case). This is possible because locally, over an affine $S=Spec(R)$
any section of $\ul{H^1}$ is given by finitely many data, so the situation
is pulled back from the spectrum of a finitely generated subring. 
For the same reason, it suffices to establish
 the additivity and normalization properties (b) and (c) 
 in the case of Noetherian $S$, in which case they are  basic properties
 of the residue symbol, formulated in \cite{Conrad} (A.1.5) and proved
 there afterwards. 
 
 \vskip .3cm
 
\noindent \underbar{Step 2: the even part:} Let $S$ be an arbitrary super-scheme.
Then both the source and target of the desired morphism $\on{Res}_{\Gamma/S}$
are $\ZZ/2$-graded, so we need to construct the even compoment
$$\on{Res}_{\Gamma/S, \even}:  \,\,q_*\bigl(\Kc_\Gamma\otimes\Omega^1_{S\times C/S}\bigr)_\even
\lra \Oc_{S, \even},$$
as well as the odd compoment $\on{Res}_{\Gamma/S, \odd}$. 
Notice that we have the ordinary scheme $\widetilde S = (\ul{S}, \Oc_{S,\even})$.
Further, since $C$ is a purely even curve, the $\ZZ/2$-grading in
$\Omega^1_{S\times C/S}$ is induced by that on $\Oc_S$, which means that
$$\bigl(\Omega^1_{S\times C/S}\bigr)_\even \,\,=\,\,\Omega^1_{\widetilde S\times C/
\widetilde S}.$$
So we define
$$\on{Res}_{\Gamma/S, \even}\,\,=\,\,\on{Res}_{\widetilde\Gamma/\widetilde S},$$
where $\widetilde\Gamma$ is the union of the graphs of the morphisms
$\widetilde c_i: \widetilde S\to C$. Parts (b) and (c) of the lemma for
$\on{Res}_{\Gamma/S, \even}$ follow.

\vskip .3cm

\noindent \underbar{Step 3: the odd part:} We now reduce to the previous case
by using a version of the ``even rules" method of 
\cite{Deligne-Morgan} \S 1.7.
Let $\Lambda[\xi]$ be the exterior algebra
in one variable, so $\Spec\,\Lambda[\xi]=\mathbb A^{0|1}$.
 For any super-commutative algebra $R$, its odd part $R_\odd$
can be identified with a subspace of the even part $(R\otimes\Lambda[\xi])_\even$,
to be precise, with
$R_\odd\cdot\xi$, which is the same as
 the kernel of the multiplication by $\xi$ in $(R\otimes\Lambda[\xi])_\even$.
 Therefore, in order to define $\on{Res}_{\Gamma/S, \odd}$, we consider
 $S^\dagger= S\times\mathbb A^{0|1}$ and morphisms $c_i^\dagger: S^\dagger\to S\to C$,
 with the union of their graphs being the super-scheme 
 $\Gamma^\dagger=\Gamma\times\mathbb A^{0|1}$. 
 We then define $\on{Res}_{\Gamma/S, \odd}$ to be the restriction of
 $\on{Res}_{\Gamma^\dagger/S^\dagger, \even}$ on the kernel of the multiplication with $\xi$
 in its source and target. Parts (b) and (c) of the lemma for
$\on{Res}_{\Gamma/S, \odd}$ follow from their validity for 
$\on{Res}_{\Gamma^\dagger/S^\dagger, \even}$. 

It remains to show the compatibility of
$\on{Res}_{\Gamma/S}$ defined in terms of its even and odd components, with
arbitrary base change of super-schemes $S'\to S$. It is enough to assume that
$S=\Spec(R)$, $S'=\Spec(R')$, so we have a morphism of super-commutative
algebras $R\to R'$. For $\on{Res}_{\Gamma/S, \even}$ this follows
from the compatibility of the Grothendieck duality with the  base change
for $R_\even\to R'_\even$, while 
 for $\on{Res}_{\Gamma/S, \odd}$ it follows from compatibility with the base change
 for $R[\xi]\to R'[\xi]$. This finishes the proof of Lemma 
 \ref{Res-Gamma/S-lemma}
  and of Proposition 
\ref{transgression-additive-proposition}.


\subsection {Additive functions on $\Lc X$ and the Radon transform.}

We start with several versions of the Poincar\'e lemma.

\begin{lem}\label{poincare-lemma-one}
 Let $Y$ be a smooth algebraic super-variety, $Z\subset Y$ be a smooth
sub-super-variety with sheaf of ideals $I_Z\subset \Oc_Y$, and 
$$\widehat{\Omega}^\bullet_Y \,\,=\,\,\pro_{m\geqslant 0} \Omega^\bullet_Y/I_Z^{m+1}
\Omega^\bullet_Y$$
be the completion of the de Rham complex of $Y$ along $Z$. Then the complex
$\widehat\Omega^\bullet_{Y|Z}=\Ker\bigl\{\widehat\Omega^\bullet_Y\to\Omega^\bullet_Z\bigr\}$
is exact everywhere on each affine open set of $Z$.
\end{lem}

\noindent {\sl Proof:}  Denote by $N^*=I_Z/I_Z^2$ the conormal bundle of $Z$ in $Y$.
Filtering by powers of $I_Z$, we equip $\widehat\Omega^\bullet_{Y|Z}$
with a decreasing complete filtration whose quotients are nothing but the homogeneous
pieces of the Koszul complex:
$$S^pN^*\lra S^{p-1}N^*\otimes_{\Oc_Z} \Lambda^1 N^*\lra S^{p-2}N^*\otimes_{\Oc_Z}\Lambda^2 N^*
\lra ..., \quad \quad p\geqslant 1.$$
Each such quotient is exact on each affine open set. \qed

\begin{lem}\label{poincare-lemma-two}
Let $X$ be a smooth super-manifold.
The relative De Rham complex
$$\pi_*\Omega^\bullet_{\Lc X|\Lc^0X}=\bigl\{
\pi_*\Oc_{\Lc X|\Lc^0X}\buildrel d\over\lra
\pi_*\Omega^1_{\Lc X|\Lc^0X}\buildrel d\over\lra
\dots\bigr\}$$
is exact everywhere on the Zariski topology of $X$.

\end{lem}

\noindent {\sl Proof:} The statement being local, we can assume that $X$ admits an
\'etale coordinate system $\phi: X\to\mathbb A^{d_1|d_2}$. We then have
a realization of $\Lc X$ as a double ind-pro-limit of the schemes $\Lc_n^\epsilon(\phi)$,
as in Remark \ref{L-n-epsilon-remark}. Fixing $m>0$, let
$$\Lc_n^m(X) \,\,=\,\,``\hskip -.5cm\ind_{\epsilon_i=0, \, i<-m}\hskip -.5cm "
\hskip .5cm\,\Lc_n^\epsilon X,$$
where the limit is taken over those $\epsilon\in\Eb$ which have $\epsilon_i=0$
for $i<-m$. Then $\Lc_n^m(X)$ is isomorphic to the
 formal neighborhood of the smooth super-algebraic
variety $\Lc^0_n X$ inside the product of $\Lc^0_n X$ with an affine super-space
of  dimension $d_1m|d_2m$. So $\Omega^\bullet_{\Lc X|\Lc^0X}$
is a complex of the kind considered in Lemma \ref{poincare-lemma-one}
and therefore it is exact on each affine open set. Now, 
$$\pi_*\Omega^\bullet_{\Lc X|\Lc^0X}\,\,=\,\,\pro_m \ind_n \,\,\pi_{n*}\Omega^\bullet_{\Lc^m_n X|\Lc^0_nX},$$
where $\pi_n: \Lc^0_n X\to X$ is the projection. Further,
the ind-pro-system has the maps in the ind-direction injective and the maps
in the pro-direction surjective. So the double limit is exact as well.
\qed

\vskip .3cm

For  a closed 2-form $\omega$ on $X$ we have a closed 1-form
$\tau(\omega)$ in $\Omega^{1, \cl}_{\Lc X|\Lc^0 X}$.  
Let $d^{-1}(\tau(\omega))$ be its unique preimage under the de Rham differential
which lies in $\Oc_{\Lc X|\Lc^0X}$. 

\begin{thm}\label{additive-2-forms-theorem}
The correspondence $\omega\mapsto d^{-1}\tau(\omega)$ defines a morphism
of sheaves $d^{-1}\tau: \Omega^{2, \cl}_X \to \Ac dd_X$,
which is an isomorphism.
\end{thm}

This theorem was proved in \cite{KV3} when $X$ is an even manifold
using the results of \cite{GMS1}.
Here we give an independent proof in the more general context
of super-manifolds.
The morphism $d^{-1}\tau$ can be called the {\em Radon transform}
on the space of formal loops. If $\omega$ is a symplectic form
on $X$, the function $d^{-1}\tau(\omega)$ is the formal loop space version of
the {\em  symplectic action functional}. 

\vskip .2cm

To prove Theorem \ref{additive-2-forms-theorem}, we associate 
to any additive function $f$ on $\Lc X$
a 2-form  as follows.
Consider the embedding of constant loops
$$\epsilon:X\hookrightarrow \Lc^0X\hookrightarrow\Lc X.$$ 
We will study the behaviour of $f$ on the first and second infinitesimal
neighbourhoods of $X$ in $\Lc X$. 
First, let us introduce the following notation
$$\Omega^1_{\Lc X}|_X=
\epsilon^{-1}(\Omega^1_{\Lc X})\otimes_{\epsilon^{-1}(\Oc_{\Lc X})}\Oc_X.$$
For a section $\omega$ of  $\Omega^1_{\Lc X}$ we denote by
$\omega|_X$ its image in $\Omega^1_{\Lc X}|_X$ and call it
the {\em restriction} of $\omega$ to $X$

\begin{lem}
(a) We have
$\Omega^1_{\Lc X}|_X=
\Omega^1_X\lb t^{-1}\rb.$

(b) Dually, defining $\Theta_{\Lc X}|_X=\Der(\Oc_{\Lc X},\Oc_X)$, 
the sheaf of continuous derivations, we have
$\Theta_{\Lc X}|_X=\Theta_X\lb t\rb$.
\end{lem}

\noindent{\sl Proof :}
Part (a). Let $f$ be a local section of $\Oc_X$. Then, for any $m\in\ZZ$,
we have that $f[m]$ is a local section of $\Oc_{\Lc X}$, and so
$d(f[m])$ is a local section of $\Omega^1_{\Lc X}$. Our identification maps
$f[m]\on{d}(g[n])$ to $(f\on{d} g)t^{m+n}$.

Part (b). Let $\xi$ be a local section of $\Theta_X$.
We denote by $\partial_\xi$ the corresponding derivation of $\Oc_X$.
Let's now define $\partial_{\xi[n]}$ to be the derivation $\Oc_{\Lc X}\to\Oc_X$
given by
$$\partial_{\xi[n]}(f[m])=\delta_{m,n}(\partial_\xi f).$$
This define a subsheaf $\Theta_X[n]$ of $\Theta_{\Lc X}|_X$.
Our identification maps $\Theta_X[n]$ to $\Theta_X\,t^n$. \qed

\vskip 3mm

The group $\GG_m\subset K$ acts on $\Lc X$ by {\it the rotation of the loop}
$t\mapsto \lambda t$.
So it acts also on the pro-sheaves
$\Omega^1_{\Lc X}|_X$ and $\Theta_{\Lc X}|_X$.
The homogeneous components of degree $n$ are respectively
$\Omega^1_X[n]=\Omega^1_X\,t^n$ and $\Theta_X[n]=\Theta_X\,t^n$.

\vskip 3mm

\begin{lem}
If $\omega\in\Ac dd^1_X$ is an additive 1-form on $\Lc X$ then
the restriction $\omega|_X$ is equal to 0. In particular, 
  if $f$ is an additive function on $\Lc X$ then
the differential $d_x f$ vanishes along $x$.
\end{lem}

\noindent {\sl Proof:} 
It is enough to prove the first claim.
Since $\omega$ is additive, it is, in particular, $(\gen, K)$-invariant.
Thus $\omega$ is invariant under the subgroup $\GG_m$, and
so is $\omega|_X \in \Omega_X\lb t^{-1}\rb$. Since
$\GG_m$ acts on $\Omega_x\cdot t^n$ via the character
$\lambda\mapsto\lambda^n$, we conclude that $\omega|_X$ should lie
in the subspace $\Omega_X\cdot t^0$. But the $t^0$-component should also
vanish, since 
the condition that $\omega=0$ on 
$$\Theta_{\Lc^0 X}|_X=\prod_{n\geqslant 0}\Theta_X[n],$$ 
is also included in
the property of being additive.\qed

\vskip .2cm

We now continue our argument. As
 the value and the differential of $f$ vanishes identically along $X$,
we have the  invariantly defined Hessian, which is a quadratic form on the 
restriction of the tangent bundle to $X$:
$$H(f): \,\, S^2 \Theta_{\Lc X|_X} \,\,\to\,\,\Oc_X.$$
Let $B(f)$ be the corresponding symmetric bilinear form. 
From the $\GG_m$-invariance of $f$ and thus of $B(f)$ we conclude that the
only possibly non-trivial homogeneous components of $B(f)$ are the pairings
$$B^n(f): \,\, \Theta_{X}[-n] \,\otimes \,\Theta_{X}[n] \,\,\to \,\,\Oc_X, \quad
n\neq 0.$$
By identifying each $\Theta_{X}[n]$ with $\Theta_{X}$, we can associate
to $B^n(f)$ a contravariant 2-tensor $\omega^n\in H^0(X, \Omega^1_X\otimes\Omega^1_X)$:
\be 
\omega^n (v, w) \,\, = \,\, B^1(f) \bigl(v[-n], w[n]\bigr). 
\ee
Here $v,w$ are vector fields on $X$. 
For example, let $X=\mathbb A^{d_1|d_2}$, so tangent vectors $v, w$ to $X$ at any point $x$
with  values in a super-commutative algebra $R$
can be seen as elements of $R\otimes \mathbb C^{d_1|d_2}$. Then
\be 
\omega_x^n(v, w) \,\,\, = \,\,\, {1\over 2} 
{d^2\over d\eps^2}\biggl|_{\eps =0} \,\,
f[\eps vt^{-n} + x + \eps wt^n]\,\,\in\,\, R. 
\ee

\begin{prop}
For each $n>0$ we have $\omega^{\pm n} = n\omega^{\pm 1}$.
\end{prop}

\noindent {\sl Proof:}
Consider the morphism
$$\Psi^n: \, \Lc X \,\to\,\Lc X$$
induced by the change of variable $t=u^n$ in the formal series. More precisely,
this change of variable induces, for any super-scheme $S$, a morphism of
super-spaces
$$\bigl(\ul{S}, \Oc_S\lb t\rb^\sqr\bigr) \,\,\, \lra \,\,\,
 \bigl(\ul{S}, \Oc_S\lb u\rb^\sqr\bigr),$$
and thus we have an endomorphism of the functor representing $\Lc X$. 
It is clear that $\Psi^n$ is identical on $X$ and its differential
has the following form on $\Theta_X[1]$:
$$d\Psi^n: \,\,\Theta_{X}[1] \,\,\,\lra\,\,\,\Theta_{X}[n], \quad v[1]\mapsto v[n].$$
So our statement would follow from the next lemma.

\begin{lem}
Any additive $m$-form $\omega\in \Ac dd_X^m$ satisfies 
$(\Psi^n)^*(\omega) = n\omega$.
\end{lem}

\noindent {\sl Proof:} Our change of variable gives a morphism
$$D\,\, =  \,\,\Spec(\CC[u])\,\lra \,
C \,\, = \,\,\Spec(\CC[t]).$$
Let $\ZZ_n$ be the cyclic group of order $n$ acting on $D^n$
and thus on $\Lc_{D^n} X$ by cyclic permutations. The morphism
$\Psi^n:\Lc X\to\Lc X$ extends to a morphism of global loop spaces
$$\widetilde{\Psi}^n: \,\, \Lc_C X \,\,\, \to\,\,\, (\Lc_{D^n}X)/\ZZ_n,$$
as the pre-image of a non-zero point $t_0\in C$ consists on $n$ points
defined up to a cyclic permutation.

\vskip .2cm

Now, the $m$-form $\omega_{D^n}\in \Omega^m_{\Lc_{D^n}X/D^n}$
is invariant under all permutations, in particular, under $\ZZ_n$
and so descends to a $m$-form
$\widetilde\omega_{D^n}$ on $(\Lc_{D^n}X)/\ZZ_n$.
Consider the $m$-form $(\widetilde\Psi^n)^*(\widetilde\omega_{D^n})$ on
the ind-scheme $\Lc_C X$.
The fiber of $\Lc_CX$ over each $t\in C$ is identified with $\Lc X$
canonically up to the action of the group scheme $K$. Now, for $t\neq 0$
the restriction of 
the $m$-form $(\widetilde\Psi^n)^*(\widetilde\omega_{D^n})$ 
to this fiber is equal to $n\omega$ because each of the $n$ pre-images
of $t$ in $D$ will contribute a summand equal to $\omega$, 
in virtue of the additivity of $\omega$.
On the other hand, for $t=0$ the restriction is equal to 
$(\Psi^n)^*(\omega)$
by definition. This proves the lemma.
\qed

\begin{prop}\label{skew-symmetry-omega-prop}
Let $f\in\Ac dd_X$ and let $\omega^1$ be defined as above. Then:

(a) The tensor $\omega^1$ is skew symmetric in the super sense,
yielding a differential 2-form on $X$. 

(b) We have $\omega^n = n\omega^1$
for all $n\neq 0$,

(c) The 2-form $\omega^1$ is closed: $d\omega^1 = 0$. 
\end{prop}

\noindent{\sl Proof :}
It is enough to assume that $X$ is affine
and is equipped with an \'etale morphism $\phi: X\to \mathbb A^{d_1|d_2}$.
We denote by $x_1, ..., x_N\in\Oc(X)$, $N=d_1+d_2$, the pullbacks under $\phi$
of the (odd and even) coordinate functions on $\mathbb A^{d_1|d_2}$. 
Then $dx_1, ..., dx_N$
form an $\Oc_X$-basis of $\Omega^1_X$, and we denote by
$\partial/\partial x_1, ..., \partial/\partial x_N$ the dual basis of
$\Theta_X$. We can then use Taylor expansions of functions on $X$ and
$\Lc(X)$ in the same way as if $X$ was a Zariski open subset in $\mathbb A^{d_1|d_2}$.
The identification
$$\Lc(X) \,\,\,\simeq \,\,\, \Lc(\mathbb A^{d_1|d_2})\times_{\mathbb A^{d_1|d_2}}X,$$
see Remark \ref{L-n-epsilon-remark}
 and \cite{KV1}, Proposition 1.6.1, means that we have the functions
$x_{i,n}$ on $\Lc(X)$, with $i=1, ..., N$ and $n\in \ZZ$ which we can think
as the coefficients of $N$ indeterminate Laurent series
\be 
x_i(t) \,\,\, = \,\,\, \sum_{n=-m}^\infty x_{i,n} t^n, \quad i=1, ..., N.
\ee
Thus we can expand the function 
$f$ in the pro-algebra $\Oc(\Lc X)$
near each $\CC$-point of  $X\subset \Lc (X)$ as a series in these coordinates.
 
Let us consider only Laurent series starting with terms with $t^{-1}$ and write
$x_n= (x_{1,n}, ..., x_{N, n})$ for the vector of the $n$th coefficients.
Then we can write
\be\label{general-expansion-f}
\begin{aligned}
f[x_{-1}t^{-1} + x_0 + x_1 t + x_2 t^2 + ...] \,\,\, = \,\,\,
\sum_{i,j} \omega_{ij}(x_0) x_{i,-1} x_{j,1} \,\, +\\
+  \sum_{(i\leqslant j), \, k} \psi_{ijk}(x_0)
x_{i,-1}x_{j,-1} x_{k,2} \,\, + \,\, 
\sum_{(i\leqslant j),\,(k\leqslant l)} \phi_{ijkl}(x_0)x_{i,-1}x_{j,-1}x_{k,1}
x_{l,1} + ...
\end{aligned}
\ee
Note that the above expansion is quasihomogeneous of degree 0 in the $x_{i,n}$,
because of the $\GG_m$-invariance of $f$. 
We identify the coordinate $x_{i,0}$ on $\Lc(X)$
with the coordinate $x_i$ on $X$.
Now, the  $\omega_{ij}(x_0)$ are  
nothing but the coefficients of the tensor $\omega^1$.
More precisely, we set
$$\omega^1 = \sum_{i,j} \omega_{ij} \on{d}x_i \otimes \on{d}x_j,
\quad\omega^1_{ij}=(-1)^{d_i}\omega_{ij}(x_0).$$ 
So our first task is to prove
the antisymmetry of the $\|\omega_{ij}\|$ in the super sense, i.e., that
$$\omega_{ij}=(-1)^{(1+d_i)(1+d_j)}\omega_{ji},\quad
d_i=\deg(x_i)\in\ZZ/2.$$
We now explain a method allowing us to exploit the additivity of $f$ in order to
obtain information about the coefficients such as $\omega_{ij}(x)$. 
Fix a $\CC$- point $o\in X$ and assume that the functions $x_i$ vanish at $o$,
so we think of $o$ as the origin of coordinates and study the behavior of $f$
near $o\in X\subset \Lc(X)$.
Let  $a=(a_1, ..., a_N)$ and $b=(b_1, ..., b_N)$ be two vectors of independent
variables of the same parities as $(x_1,\dots,x_N)$
which we eventually suppose to be nilpotent of some degree $d$,
so we define
$$R \,\,\, = \,\,\,
 \CC\bigl[a_i, b_i|\, i=1, ..., N\bigl]/ \bigl(a_i^d, b_i^d|\, i=1, ..., N\bigl).$$
Consider the  rational loop
\be\label{gamma-t-equation}
\gamma(t) \,\,\, = \,\,\,{a\over t} + {b\over\lambda-t},
\ee
where $\lambda\in\CC$ is a parameter. To be precise, $\gamma(t)$
is the unique $R[\lambda]$-point of
$\Lc_{\mathbb A^2}X$ whose image under $\phi: X\to\mathbb A^{d_1|d_2}$
is the rational loop in the right hand side of \eqref{gamma-t-equation}.
Note that the canonical map
$\Lc_{\mathbb A^2}X\to\mathbb A^2$ takes $\g(t)$ to the $R[\lambda]$-point $(0,\lambda)$ of
 $\mathbb A^2$.
Since $f$ is additive,  we have the function $f_{\mathbb A^2}$, whose value at $\gamma(t)$ is an 
element of $R[\lambda]$. 
On the other hand, we can expand $\gamma(t)$ at each of the two poles,
which gives:
$$\gamma(t) \,\,\, = \,\,\, at^{-1} + {b\over\lambda} t^0 + {b\over\lambda^2}t + {b\over\lambda^3} t^2 + ...$$
near $t=0$. Now, near $t=\lambda$ we have the coordinate $s=\lambda-t$, and
$$\gamma(t) \,\,\, = \,\,\, bs^{-1} + {a\over\lambda} s^0 + {a\over\lambda^2}s^1 + {a\over\lambda^3} s^2 + ...$$
We see that the coefficients of each individual expansion become singular as $\lambda\to 0$, but the
value
$$f_{\mathbb A^2}[\gamma(t)]\,\,\, =\,\,\,  f\biggl[at^{-1} + {b\over\lambda} t^0 +
{b\over\lambda^2}t + {b\over\lambda^3} t^2 + ...\biggr]
\,\,+ \,\, f\biggl[bs^{-1} + {a\over\lambda} s^0 + {a\over\lambda^2}s^1 + {a\over\lambda^3} s^2 + ...\biggr]$$
must be regular at $\lambda =0$. 
So expanding each summand into a Taylor series using
\eqref{general-expansion-f},  we have that the coefficients at 
each negative power of $\lambda$ must cancel, which
provides a system of constraints on the 
coefficients $\omega_{ij}$, $\psi_{ijk}$ etc.
Thus, we have:
$$f\biggl[at^{-1} + {b\over\lambda} t^0 +
{b\over\lambda^2}t + {b\over\lambda^3} t^2 + ...\biggr] \,\,\, = \,\,\,
\sum_{i,j} \omega_{ij}\biggl({a\over\lambda}\biggr) {a_ib_j\over\lambda^2}  \,\, + \,\,
\sum_{(i\leqslant j),\,k} \psi_{ijk}\biggl({b\over\lambda}\biggr) {a_ia_jb_k\over\lambda^3} + \cdots$$
where dots stand for terms with $1/\lambda^4$ and higher. 
To arrive at the precise 
coefficients at the powers of $\lambda$, 
we need to further  expand $\omega_{ij}(b/\lambda)$,
$\psi_{ijk}(b/\lambda)$ etc. near the point $o$, 
using the Taylor formula, which gives:
$$\gathered
f\biggl[at^{-1} + {b\over\lambda} t^0 +
{b\over\lambda^2}t + {b\over\lambda^3} t^2 + ...\biggr] \,\,\, = \,\,\,
\sum_{i,j} \omega_{ij}(o) {a_ib_j\over\lambda^2} \,\, + \cr
+\,\,\sum_{i, j,k} {\partial\omega_{ij}
\over\partial x_k}(o) {b_ka_ib_j\over\lambda^3} \,\, + \sum_{(i\leqslant j),\, k} \psi_{ijk}(o) {a_ia_jb_k\over
\lambda^3} + \cdots,
\endgathered$$
and similarly for the other summand.
So the cancellation of the terms with $1/\lambda^2$ 
in $f_{\mathbb A^2}[\gamma(t)]$ implies
that 
$$\omega_{ij}(o) + (-1)^{\deg(a_i)\,\deg(b_j)}\omega_{ji}(o)=0. $$
Since $d_i=\deg(a_i)$, $d_j=\deg(b_j)$ and $o$
can be any point of $X$,
this proves the antisymmetry of $\omega^1$ 
and thus parts (a) and (b) of Proposition \ref{skew-symmetry-omega-prop}. 

\vskip .2cm

Continuing further, for $j\leqslant k$,  cancellation of 
the coefficients at $a_ib_jb_k/\lambda^3$ gives
\be\label{omega-psi-equation}
\gathered
(-1)^{d_id_k+d_jd_k}
{\partial 
\omega_{ij}\over\partial x_k}(o) + 
(-1)^{d_id_j}
{\partial \omega_{ik}\over\partial x_j}(o)
+(-1)^{d_id_j+d_id_k}\psi_{jki}(o)=0.
\endgathered
\ee
So the terms with the derivatives of $\omega_{ij}$ 
become mixed with the terms with $\psi_{ijk}$.
To avoid this mixing, we modify our approach by 
considering the rational loop with three poles
\be 
\delta(t) \,\,\, = \,\,\, {a\over t} + {b\over \lambda-t} + 
{c\over \lambda +t}, 
\ee
where $c=(c_1, ..., c_N)$ is a third group of nilpotent independent variables
of the same parities as $(x_1,\dots,x_N)$. 
As before, $f_{\mathbb A^3}[\delta(t)]$ is the sum of values of $f$ at the three expansions of
$\delta(t)$: near $t=0$ where it is
$$\delta(t) \,\,\, = \,\,\, at^{-1} + {b+c\over\lambda}t^0  + {b-c\over\lambda^2} t + {b+c\over\lambda^3}
t^2 + ...,$$
near $t=\lambda$, where the expansion in $s=\lambda-t$ is
$$\delta(t) \,\,\, = \,\,\, bs^{-1} + {a+c/2\over \lambda}s^0 + {a+c/4\over\lambda^2} s
+ {a+c/8\over\lambda^3}s^2 + ...,$$
and near $t=-\lambda$, where the expansion in $u= t+\lambda$ is
$$\delta(t) \,\,\, = \,\,\, cu^{-1} + {-a+b/2\over\lambda} s^0 + {-a+b/4\over\lambda^2} s + {-a+b/8\over
\lambda^3} u^2 + ...$$
The sum of the values of $f$ at these three expansions should not have terms with negative powers of
$\lambda$. As before, the cancellation of the terms with $1/\lambda^2$ gives the antisymmetry of the $\omega_{ij}$,
while the coefficient at $1/\lambda^3$ is found by using the Taylor formula to be:
$$\sum_{i,j} {\partial \omega_{ij}\over\partial x_k}(o) (b_k+c_k) a_i (b_j-c_j) \,\, + \,\,
\sum_{(i\leqslant j), k} \psi_{ijk}(o)a_ia_j(b_k-c_k) \,\, +$$
$$+\sum_{i,j,k} {\partial \omega_{ij}\over\partial x_k}(o)(a_k + c_k/2)b_i (a_j+c_j/4) \,\, + \,\,
\sum_{(i\leqslant j), k} \psi_{ijk}(o) b_i b_j (a_k-c_k/4) \,\,+$$
$$+ \sum_{i,j,k} {\partial \omega_{ij}\over\partial x_k}(o)(-a_k+ b_k/2)c_i(-a_j+b_j/4) \,\,+ \,\,
\sum_{(i\leqslant j), k} \psi_{ijk}(o) c_ic_j(-a_k+b_k/4).$$
In this sum we concentrate on the mixed monomials of the form $a_i b_j c_k$.
The coefficient at such a monomial is found to be
$${1\over 2}\biggl( (-1)^{d_i+d_id_k+d_jd_k}{\partial \omega_{ij}^1\over\partial x_k}(o)+
(-1)^{1+d_i+d_id_j}{\partial \omega_{ik}^1\over\partial x_j}(o)+
(-1)^{d_j}{\partial\omega_{jk}^1\over\partial x_i}(o)\biggr).$$
So vanishing of such coefficients implies that $\omega^1$ is closed, because
$$\gathered
\on{d}\omega^1=\sum_{i,j,k}\biggl( (-1)^{(1+d_k)(d_i+d_j)}{\partial \omega_{ij}^1\over\partial x_k}(o)+
(-1)^{(1+d_i)(1+d_j)}{\partial \omega_{ik}^1\over\partial x_j}(o)+\cr
\hskip2cm
+
{\partial\omega_{jk}^1\over\partial x_i}(o)\biggr)\on{d}x_i\on{d}x_j\on{d}x_k.
\endgathered$$
Proposition \ref{skew-symmetry-omega-prop} is proved. 
\qed

\vskip .3cm

We will denote the 2-form $\omega^1$ simply by $\omega$ and call it
the {\em tangential 2-form} of $f$. To emphasize its dependence on $f$,
we will write $\omega = Df$.

\begin{lem}
Let $\omega'\in\Omega^{2, \cl}(X)$ be a
given closed 2-form
and let $f = d^{-1}\tau(\omega')$. Then $Df = \omega'$. 
\end{lem}

\noindent {\sl Proof:} As earlier, it is enough to prove the statement
in the formal neighborhood of any point $x\in X$, and so
the statement reduces to that for the formal completion of $\mathbb A^{d_1|d_2}$ at 0.
Because of the formal Poincar\'e lemma, we can assume that $\omega' = d\eta$
is exact, so $f = \tau(\eta)$. 
Let $x_1, ... , x_N$ be the (even and odd) coordinates in $\mathbb A^{d_1|d_2}$, so we write
$\eta=  \sum_{i=1}^N \eta_i(x) dx_i$. Here $\eta_i(x)\in \CC\lsb x_1, ..., x_N\rsb$
is a formal power series in the even variables with coefficients being  elements
 of the exterior algebra
in the odd variables.
Then
$$\omega' \,\,\, = \,\,\, \sum_{i<j} 
\omega'_{ij}(x) dx_i dx_j, \quad \omega'_{ij}(x)
\,\,\,  = \,\,\,
{\partial \eta_i\over\partial x_j} - {\partial \eta_j\over\partial x_i}.$$
Let $\omega = Df$, so we need to prove that $\omega = \omega'$. 
Let $e_1, ..., e_N$ be the basis of $\CC^{d_1|d_2}$ corresponding to the
coordinate system $x_1, ..., x_N$. We view $\CC^{d_1|d_2}$ as the tangent space to
$\mathbb A^N$ at any $\CC$-point. Then we need to prove that
$$\omega'_{ij}(x) \,\,\, = \,\,\, \omega_x(e_i, e_j) \,\,\, := \,\,\,{1\over 2} 
{d^2\over d\eps^2}\biggl|_{\eps =0} f[\eps e_i t^{-1} + x + \eps e_jt].$$
For the purposes of such a proof all the coordinates $x_k, \, k\neq i,j$
appear as parameters (constants with respect to the differentiation),
so we can assume that $N=2$, $i=1, \,\, j=2$. By splitting $\eta$ into two
summands and switching the roles of $x_1$ and $x_2$, it is enough to assume
that $\eta = \eta_2(x_1, x_2) dx_2$. Further, by decomposing $\eta_2$
into monomials, we reduce to the case  $\eta = x_1^a x_2^b dx_2$. 
The formal loop 
$$\gamma(t) \,\, = \,\,\gamma_\eps (t) \,\, = \,\, \eps e_1 t^{-1} + x + \eps e_2 t$$
has the coordinates
$$x_1 [\gamma(t)]\,\, = \,\,  \eps t^{-1} + x_1, \quad 
x_2[\gamma(t)] \,\, = \,\, x_2 + \eps t.$$
So we have 
$$f[\gamma(t)] \,\, = \,\, \Res_{t=0} \,\, \biggl[(\eps t^{-1} + x_1)^a (x_2+\eps t)^b d (x_2+\eps t)
\biggr]
\,\, = \,\, \eps^2 a x_1^{a-1} x_2^b,$$
see \cite{KV3}, Example (1.3.8).
This is exactly $\eps^2$ times
the coefficient at $dx_1 dx_2$ of 
$$d(x_1^a x_2^b dx_2)=d\eta=\omega'.$$ 
\qed

\vskip .3cm

We have proved that the morphism of sheaves
$$D: \,\, \Ac dd_X \,\,\lra \,\, \Omega^{2,\cl}_X$$
is left inverse to $d^{-1}\tau$, 
so $D$ is surjective and $d^{-1}\tau$ is injective.
To prove that they are mutually inverse isomorphisms, it suffices to prove the following.

\begin{prop}
If an additive function $f$ is such that $Df=0$ identically,
then $f=0$ identically.
\end{prop}

\noindent {\sl Proof, first step:} we prove that $f[x(t)]=0$ if $x(t)$
is any formal loop whose expansion begins with terms with $t^{-1}$, 
and so is given by an expansion as in \eqref{general-expansion-f}.
Suppose that $\omega=Df$ vanishes identically.
Then the first group of terms in the right hand side of \eqref{general-expansion-f} vanishes.
We prove inductively that all the coefficients in this expansion vanish, 
using the vanishing of the coefficient at each negative power of
$\lambda$ in $f[\gamma(t)]$. Indeed, identical vanishing of each
$\omega_{ij}$ implies, by \eqref{omega-psi-equation}, that each
$\psi_{ijk}(o)= 0$. Here $o$ can be any point, so $\psi_{ijk}\equiv 0$.
Next,  comparing coefficients at $1/\lambda^4$
in $f[\gamma(t)]$, we get a relation between the values of $\phi_{ijkl}$, the
first derivatives of  the $\psi_{ijk}$ and the second derivatives of the $\omega_{ij}$
at any given point $o$. This implies that each $\phi_{ijkl}\equiv 0$, and so on.

\vskip .3cm

\noindent {\sl Second step:} Any formal loop, i.e., each $R$-point of $\Lc(X)$
$$x(t) \,\,\, = \,\,\,\sum_{n=-M}^\infty x_nt^n,\quad x_n = (x_{1,n}, ..., x_{N,n}), \quad
x_{i,n}\in R,$$
with order of pole $M\geqslant 2$, can be deformed into a 1-parameter family 
of rational loops
each having $M$ poles of first order, by considering the $R[\lambda]$-point 
of $\Lc_{\mathbb A^M}(X)$ given by
$$x_\lambda(t) = \sum_{p=2}^M {x_{-p}\over t(t+\lambda) ... (t+(p-1)\lambda)} \,\, 
+\sum_{n=-1}^\infty x_n t^n.$$
Then the value $f_{\mathbb A^M}[x_\lambda(t)]\in R[\lambda]$ vanishes for $\lambda\neq 0$,
since then $x_\lambda(t)$ has only first order poles. Therefore
the specialization of $f_{\mathbb A^M}[x_\lambda(t)]$ to $\lambda =0$, i.e.,
$f[x(t)]$, vanishes as well.\qed

This finishes the proof of Theorem \ref{additive-2-forms-theorem}.


\subsection{ Additive forms on $\Lc X$.}

The goal of this subsection is to prove the following theorem.

\begin{thm}\label{additive-forms-transgression-theorem}
The morphism of complexes 
$\tau_X^{\geqslant 2}: \Omega^{\geqslant 2}_X\to\Ac dd^{\geqslant 1}_X[-1]$
is a quasi-isomorphism.
\end{thm}

\noindent{\sl Proof :}
Since $\Omega^\bullet_{\Lc X}=\Oc_{\Sc\Lc X}=\Oc_{\Lc \Sc X}$, and similarly for
$\Lc_{C^p}X$ for any $p$,  see Remark \ref{LSX=SLX-remark},
we conclude that $\Ac dd^\bullet_X=\Ac dd^0_{\Sc X}$.
Further, the De Rham differential in $\Ac dd^\bullet_X$ is just the action
on $\Ac dd^0_{\Sc X}$ of the vector field $D$ discussed in Subsection
\ref{N=1-supersymmetry-subection}.
Next, by Theorem \ref{additive-2-forms-theorem} applied to $\Sc X$
we have a sheaf isomorphism
$$d^{-1}\tau: \Omega^{2, \cl}_{\Sc X} \lra \Ac dd^0_{\Sc X},$$
given by the transgression on $\Sc X$.
Now, by Corollary \ref{complex-corollary} (case $p=2$), we have
a derived category isomorphism
$$q: \ten_{\geqslant 3}(\Omega^\bullet_X) \lra \Omega^{2,\cl, \bullet}_{\Sc X},$$
represented by the diagram \eqref {two-quasiisomorphisms-equation}
of quasi-isomorphisms of complexes, in our case by
\be\label{two-projections-omega-equation}
\ten_{\geqslant 3}(\Omega^\bullet_X)\lla W^\bullet\lra \Omega^{2,\cl, \bullet}_{\Sc X},
\ee
read from left to right. Here $W^\bullet$ is the total complex of the double
complex $W^{\bullet\bullet}$ defined in \eqref{double-complex-W-equation}, and the arrows
are the projections to the two edges.

 We now regard the transgression as a morphism of truncated complexes
 $$\tau_X^\ten: \, \ten_{\geqslant 3} (\Omega^\bullet_X) \lra(\ten_{\geqslant 2} \Ac dd^\bullet_X)[-1]\,\,=\,\,
 \bigl\{ \Ker(d) \lra \Ac dd^1_X\buildrel d\over\longrightarrow  \Ac dd^2_X \lra ...\bigr\}.$$
 Note that $\Ker(d)$ above is identified with the sheaf of additive functions, as we discussed
 already just before the statement of Theorem \ref{additive-2-forms-theorem}. Therefore
 we will view $\tau_X^\ten$ as a morphism of complexes
 $$\tau_X^\ten: \ten_{\geqslant 3} (\Omega^\bullet_X) \lra \Ac dd^\bullet_X[-1].$$
 
 \begin{lem}\label{commutativity-transgression-lemma}
The following diagram commutes in the derived category:
$$\xymatrix{\ten_{\geqslant 3}(\Omega^\bullet_X)\ar[r]^-q\ar[d]_{\tau^\ten_X}&\Omega^{2,\on{cl},\bullet}_{\Sc X}
\ar[d]^{d^{-1}\tau_{\Sc X}}\cr
\Ac dd^{\bullet}_X[-1]\ar[r]^s&\Ac dd^{0,\bullet}_{\Sc X}[-1].}
$$
\end{lem}

The lemma implies that $\tau_X^\ten$ is a quasi-isomorphism, since $q$ is an isomorphism in
the derived category and the other two arrows in the diagram are isomorphisms of complexes. 
Further, we deduce that  $\tau_X^{\geqslant 2}$ is a quasi-isomorphism. Indeed, the only difference between
$\ten_{\geqslant 3}(\Omega^\bullet_X)$ and $\Omega_X^{\geqslant 2}$ is the lowest degree term $\Omega^{2, \cl}_X$
attached on the left. However, the two projections in \eqref{two-projections-omega-equation}
are in fact isomorphisms on this lowest degree term,  which allows us to conclude that 
 $\tau_X^\ten$ 
will still induce a quasi-isomorphism after discarding the lowest degree terms. This induced
quasi-isomorphism is $\tau_X^{\geqslant 2}$.

\vskip .3cm

\noindent {\sl Proof of the lemma:} Consider the transgression for differential forms on $\Sc X$
$$\tau_{\Sc X}: \Omega^{\bullet \bullet }_{\Sc X} \lra \Ac dd^{\bullet -1, \bullet}_{\Sc X}.$$
Here the first grading is by the degree of differential forms on $\Sc X$ or $\Lc\Sc X$, while the
second degree is induced by the $\GG_m$-action on $\Sc X$. It is clear that $\tau_{\Sc X}$
is in fact a morphism of double complexes of degree $(-1,0)$. Indeed, we saw already that
it commutes with $D_1$, the de Rham differentials on forms on $\Sc X$ and $\Sc X$. The commutativity
with $D_2$ which is the action of the homological vector field $D$ on $\Sc X$, follows
by naturality of transgression. Note that $\Ac dd^{\bullet\bullet}_{\Sc X}$ is an $N=2$ supersymmetric
complex, since it is the additive part of the double de Rham complex of $\Lc X$. 

 Now, the quasi-isomorphism $q$ is induced by the two edge
projections \eqref{two-projections-omega-equation} of the double complex $W^{\bullet\bullet}$
obtained from $\Omega^{\bullet\bullet}_{\Sc X}$ by truncating in degrees $\geqslant 2$ for the first
grading and adding $\Ker(D_1)$, as described in \eqref{double-complex-W-equation}. 
Applying $\tau_{\Sc X}$ to $W^{\bullet\bullet}$ term by term, we map it into a
similar double complex formed out of $\Ac dd^{\bullet-1, \bullet}_{\Sc X}$ by truncating in degrees
$\geqslant 1$ for the first grading and adding $\Ker(D_1)$. Denote this complex by $W_\Sc^{\bullet-1, \bullet}$
and its total complex by $W_\Sc^\bullet$.
We have the diagram of  edge projections
\be\label{W-S-projections-equation} 
\ten_{\geqslant 2}\Ac dd_{X} \lla W^\bullet_\Sc\lra \Ac dd^{0,\bullet}.
\ee
Since $\Ac dd^{\bullet\bullet}_{\Sc X}$ is an $N=2$ supersymmetric
complex, these projections are quasi-isomorphisms. Moreover, 
the morphism in the derived category obtained by reading this diagram from left to right
is the same as the isomorphism of complexes $s$. Now, to prove the commutativity of the
diagram in the lemma, involving $q$ and $s$, it suffices to note that $\tau_{\Sc X}$ gives
a morphism of the diagram \eqref{two-projections-omega-equation} defining $q$,  to the diagram
\eqref{W-S-projections-equation} defining $s$. This finishes the proof of Lemma 
\ref{commutativity-transgression-lemma} and Theorem \ref{additive-forms-transgression-theorem}.

\begin{cor}\label{omega-2-3-corollary}
The transgression defines a quasi-isomorphism of  the 2-term complexes
$$\bigl\{\Omega^2_X\to\Omega^{3,{\cl}}_X\bigr\} \buildrel \tau\over\lra
\bigl\{ \Ac dd^1_X\to\Ac dd^{2,{\cl}}_X\bigr\}.\qed$$
\end{cor}

Note that $\tau$ is not an isomorphism of complexes. For example, any
(not necessarily
antisymmetric) contravariant 2-tensor on $X$ can  be transgressed to an additive
1-form. 


 \end{document}